# THREE DISTRIBUTIONS IN THE EXTENDED OCCUPANCY PROBLEM

**BEN O'NEILL**

Research School of Population Health (Building 62),
Australian National University,
Canberra ACT 0200

ben.oneill@anu.edu.au; ben.oneill@hotmail.com

**Abstract**
The classical and extended occupancy distributions are useful for examining the number of occupied bins in problems involving random allocation of balls to bins. We examine the extended occupancy problem by framing it as a Markov chain and deriving the spectral decomposition of the transition probability matrix. We look at three distributions of interest that arise from the problem, all involving the noncentral Stirling numbers of the second kind. These distributions give a useful generalisation to the binomial and negative-binomial distributions. We examine how these distributions relate to one another, and we derive recursive properties and mixture properties that characterise the distributions.

**Keywords**
EXTENDED OCCUPANCY PROBLEM; MARKOV CHAIN; PURE BIRTH PROCESS; OCCUPANCY DISTRIBUTION; NEGATIVE-OCCUPANCY DISTRIBUTION; SPILLAGE DISTRIBUTION.

## 1. Introduction

The classical occupancy distribution is an important and underappreciated discrete probability distribution, which describes the behaviour of the number of occupied bins when we allocate $n \in \mathbb{N}$ balls at random to $m \in \mathbb{N}$ bins. Analysis of the distribution can be found in Harkness (1969), Uppuluri and Carpenter (1971), Johnson and Kotz (1977), Kolchin, Sevast'yanov and Christyakov (1978) and Holst (1986), and some further analysis and computational aspects of are discussed in O'Neill (2021). The distribution is useful in problems involving sampling with replacement, and it is especially useful in the context of bootstrapping techniques, where it can be applied to find the probability of any given level of coverage of the original sample in a random resampling. The classical distribution covers the case where balls allocated to bins automatically "occupy" those bins. For reasons that will become clear later, it turns out to be very useful to generalise to the case where each ball has some fixed probability $0 \leq \theta \leq 1$ of "occupying" its allocated bin, and corresponding probability $1 - \theta$ of "falling through" the bin, so that it does not occupy the bin (see e.g., Uppuluri and Carpenter 1971, Samuel-Cahn 1974).

In this paper we will examine the extended occupancy problem, which seeks the marginal and conditional distributions of the "occupancy number" (counting the number of occupied bins) under conditions where the balls can fall through the bins with a fixed probability. We will derive three important distributional forms arising in this framework, and we will see how these distributions relate to one another. The first two distributions we examine are generalisations of the binomial and negative binomial distributions, and a third is a new distribution relating the extended occupancy distribution to the binomial in a useful way. All three distributions involve the noncentral Stirling numbers of the second kind, and are mathematically interesting forms that arise as norms of well-known summation formulae involving the Stirling numbers.

We examine the extended occupancy process by framing the sequence of occupancy numbers as a Markov chain, and analysing the transition matrix of the chain. To set up our analysis, we first describe the underlying mathematics of this stochastic process, based on a sequence of randomly allocated balls that can occupy or fall through the bins. Consider two independent sequences of random variables:

$$\tilde{U}_1, \tilde{U}_2, \tilde{U}_3, \ldots \sim \text{IID Unif}\{1, \ldots, m\},$$
$$Q_1, Q_2, Q_3, \ldots \sim \text{IID Bern}(\theta).$$

The first sequence represents balls allocated at random to $m$ bins and the second sequence gives indicators of whether these balls "occupy" those bins (as opposed to "falling through" the bins). From these underlying sequences we define the occupancy of each ball by the values:[1]

$$U_i = \begin{cases} \bullet & \text{if } Q_i = 0, \\ \tilde{U}_i & \text{if } Q_i = 1. \end{cases}$$

The outcome $U_i = \bullet$ means that the ball fell through its bin and so it makes no contribution to the occupancy, whereas an outcome $U_i = 1, \ldots, m$ means that the ball occupies its allocated bin.[2] For any number of balls $n$ we define the occupied **bin counts** over bins $\ell = 1, \ldots, m$ and the corresponding **occupancy number** respectively by:

$$N_{n,\ell} \equiv \sum_{i=1}^{n} \mathbb{I}(U_i = \ell) \qquad K_n \equiv \sum_{\ell=1}^{m} \mathbb{I}(N_{n,\ell} > 0)\,.$$

We have $n$ balls in our problem, but there are $n_{\text{eff}} \equiv n - N_{n,\bullet} = \sum_{i=1}^{n} \mathbb{I}(U_i \neq \bullet)$ **effective balls** (i.e., balls that occupy their allocated bins). The occupancy number counts the number of bins that are occupied, which are the bin counts that are above zero.

The occupancy process is illustrated in Figure 1 below, where we show a tabular arrangement of balls randomly allocated to bins. We show outcomes of $n = 10$ balls randomly allocated to $m = 12$ bins. Yellow squares show balls that fell through their bins and black squares show balls that occupy their bins. The bottom row of the figure shows the bin counts $N_{n,\ell}$, which add up the number of black squares in the columns above. The effective number of balls $n_{\text{eff}}$ is obtained by counting the number of black squares in the whole figure, and the occupancy number $K_n$ is obtained by counting the number of bins with at least one occupying ball (i.e., the number of columns with at least one black square).

[1] The setup for the observed data in this problem is similar to that occurring in problems with missing data. The outcome $U_i = \bullet$ is akin to the ball being "missing".

[2] Since balls that fall through their bins make no contribution to the occupancy, we can legitimately ignore loss of information of which ball they were allocated to in this case.

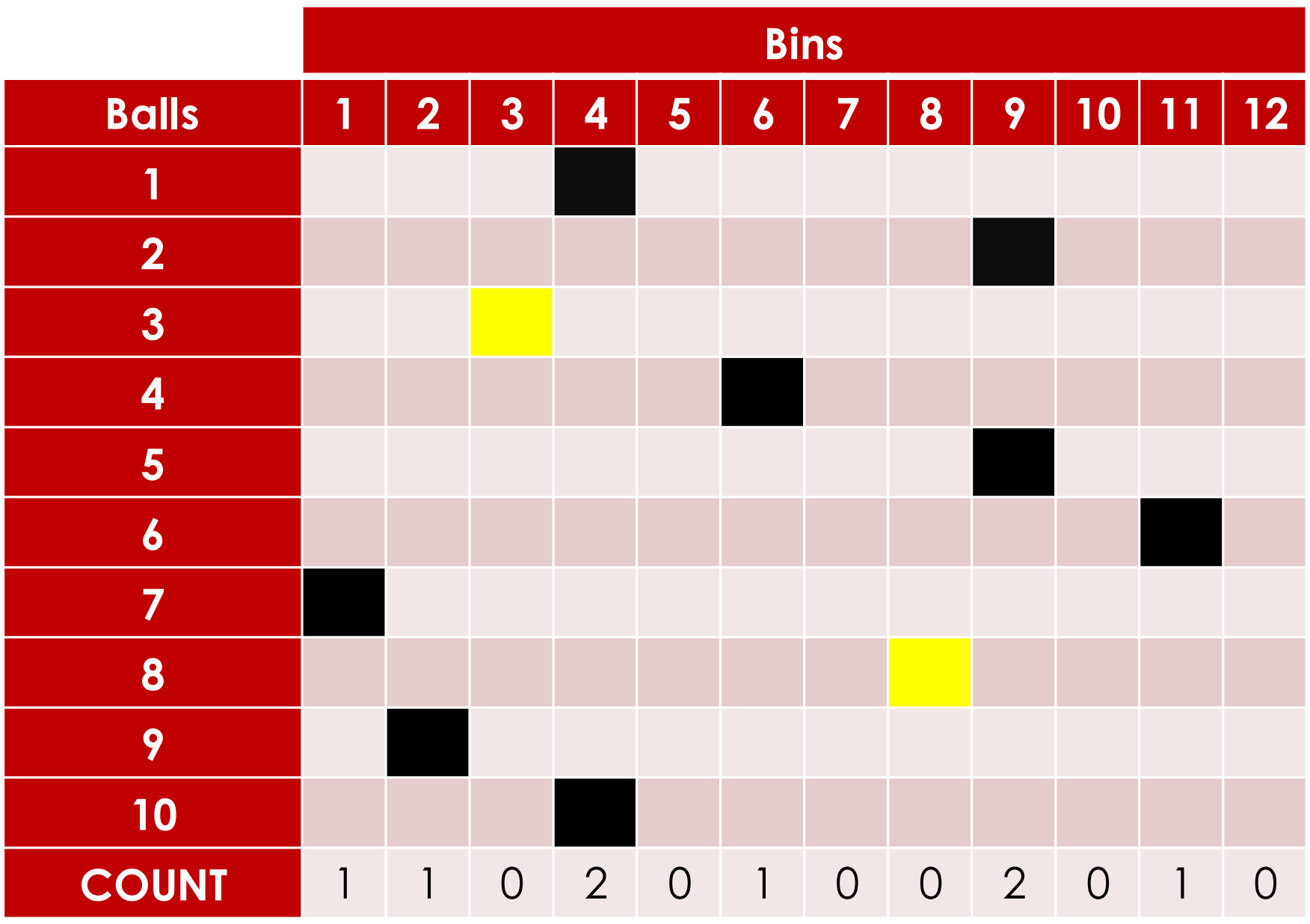

**Figure 1:** Outcomes of $n = 10$ balls randomly allocated to $m = 12$ bins. Yellow squares show balls that fell through their bins and black squares show balls that occupy their bins. Counts for each bin are shown in the bottom row. There are $n_{\text{eff}} = 8$ effective balls (black squares) in this case and the occupancy number is $K_n = 6$ (number of columns with at least one black square).

The above figure can be extended to accommodate more balls and/or bins, and the probability of a ball falling through its allocated bin can also be varied. In any case, we have now shown the mathematical foundation of the extended occupancy process and so we are in a position to state the extended occupancy problem, which seeks the marginal and conditional distributions of the occupancy number. Although we will ultimately be interested in three distributional forms arising in the extended occupancy process, our first task will be to derive a distribution form that solves this problem, and examine its properties.

**Definition (The extended occupancy problem):** For $0 \leq t \leq k \leq \min(n, m)$ we wish to find distributional forms for the marginal and conditional probabilities:

$$\mathbb{P}(K_n = k) \qquad \mathbb{P}(K_{\acute{n}+n} = k | K_{\acute{n}} = t).$$

This problem is an extension of the classical occupancy problem, which occurs when $\theta = 1$ (i.e., all balls occupy their allocated bins with probability one). □

## 2. The occupancy process and the extended occupancy distribution

Our approach to the occupancy problem is to look at the stochastic process $\{K_n | n = 0,1,2, \ldots\}$, which shows the evolution of the occupancy number as we add more balls to the process. This approach is also used in Harkness (1969) and in Uppuluri and Carpenter (1971). Each time we allocate one new ball, this ball will either occupy a bin that is not already occupied (increasing the occupancy number by one) or it will fall through its allocated bin or be allocated to a bin that is already occupied (leaving the occupancy number unchanged). Since we are assuming that balls are allocated to bins by a uniform distribution, the probability of these two outcomes is conditional on the present occupancy number, but it is not affected by which bins are occupied. In this case, the conditional probability that a newly allocated ball increases the occupancy number, given the "history" of the process, depends only on the present occupancy number, so the chain obeys the Markov property. We formalise this argument, and give the resulting transition probability and transition probability matrix, in the following theorem.

**THEOREM 1 (Markov characterisation):** Let $\boldsymbol{u}_n \equiv (u_1, u_2, \ldots, u_n)$ denote the outcomes of the first $n$ balls in the series and let $K_n = t$ denote their occupancy number. Then the conditional probability for the occupancy number with one more allocated ball is:

$$\mathbb{P}(K_{n+1} = r + t | \boldsymbol{u}_n) = \mathbb{P}(K_{n+1} = r + t | K_n = t) = \begin{cases} 1 - \theta(1 - t/m) & \text{for } r = 0, \\ \theta(1 - t/m) & \text{for } r = 1, \\ 0 & \text{otherwise.} \end{cases}$$

This conditional probability depends on the allocation history only through $K_n$, so the process satisfies the Markov property —i.e., it is a Markov chain.

**COROLLARY (Transition probability matrix):** For all the possible states $k = 0,1,2, \ldots, m$ the transition probability matrix for the Markov chain is the $(m + 1) \times (m + 1)$ bidiagonal matrix:

$$\mathbf{P} \equiv \frac{1}{m} \begin{bmatrix} m - m\theta & m\theta & 0 & \cdots & 0 & 0 \\ 0 & m - (m-1)\theta & (m-1)\theta & \cdots & 0 & 0 \\ 0 & 0 & m - (m-2)\theta & \cdots & 0 & 0 \\ \vdots & \vdots & \vdots & \ddots & \vdots & \vdots \\ 0 & 0 & 0 & \cdots & m - \theta & \theta \\ 0 & 0 & 0 & \cdots & 0 & m \end{bmatrix}.$$

(In accordance with standard conventions for Markov chains, we will index the elements of this matrix by the state values, so the first row and column will each use the zero index; i.e., the row and column indices will both run over $i, j = 0,1, \ldots, m$.)

Theorem 1 shows that the stochastic process $\{K_n | n = 0,1,2, \ldots\}$ for the occupancy number is a Markov chain with a bidiagonal transition matrix, giving us a discrete "pure birth" process. (In the nomenclature of pure birth process, an increase in the occupancy number is a "birth".) The intuition behind this transition matrix is straight-forward: if the existing occupancy number is $t$ then the newly allocated ball increases the occupancy number by one so long as it is allocated to an unoccupied bin, and does not "fall through" that bin; the probability of allocation to an unoccupied bin is $(1 - t/m)$ and the probability that the ball does not fall through the bin is $\theta$. If the newly allocated ball does falls through its allocated bin, or is allocated to a bin that is already occupied, the occupancy number does not increase.

The extended occupancy problem seeks both the marginal and conditional probabilities for the occupancy number. In the marginal problem our starting point for the chain is $K_0 = 0$, and in the conditional problem our starting point is $K_{\acute{n}} = t$. In either case, the required probabilities are easily obtained as appropriate elements of the powers of the transition matrix. Specifically, for all $k, t = 0,1,2, \ldots, m$ and all $\acute{n}, n = 0,1,2, \ldots$ we have:

$$\mathbb{P}(K_n = k) = [\mathbf{P}^n]_{0,k} \qquad \mathbb{P}(K_{\acute{n}+n} = k | K_{\acute{n}} = t) = [\mathbf{P}^n]_{t,k}.$$

Obtaining these probabilities requires us to take arbitrarily large powers of the transition matrix **P**, so it is useful to examine its spectral decomposition. It turns out that the transition matrix has eigenvectors that do not depend on the probability parameter $\theta$, and this makes our matrix characterisation especially useful, leading to a simple spectral form for the distribution.

**THEOREM 2 (Spectral decomposition):** The probability matrix **P** has **eigenvalue matrix**:

$$\mathbf{\Lambda} \equiv \begin{bmatrix} \lambda_0 & 0 & 0 & \cdots & 0 \\ 0 & \lambda_1 & 0 & \cdots & 0 \\ 0 & 0 & \lambda_2 & \cdots & 0 \\ \vdots & \vdots & \vdots & \ddots & \vdots \\ 0 & 0 & 0 & \cdots & \lambda_m \end{bmatrix} \qquad \lambda_i = 1 - \left(1 - \frac{i}{m}\right) \cdot \theta.$$

Its (unscaled) **eigenvector matrix v** and **inverse eigenvector matrix** $\mathbf{w} = \mathbf{v}^{-1}$ have elements given respectively by:

$$(\mathbf{v})_{i,j} = v_{i,j} = \binom{m-i}{j-i} \qquad (\mathbf{w})_{i,j} = w_{i,j} = (-1)^{i-j} \binom{m-i}{j-i}.$$

(Again, we remind the reader that we use indexing were the first row and column each use the zero index; this applies also to the indices for the eigenvalue and inverse eigenvalue matrices.) The columns of the eigenvector matrix (and its inverse matrix) are linearly independent, so the transition matrix is diagonalisable, with **spectral decomposition** $\mathbf{P} = \mathbf{v}\mathbf{\Lambda}\mathbf{w}$.

Theorem 2 gives the eigenvalue and eigenvector matrices of the transition matrix, which allows us to take arbitrarily large powers of this matrix using its spectral decomposition $\mathbf{P}^n = \mathbf{v}\boldsymbol{\Lambda}^n\mathbf{w}$. Uppuluri and Carpenter (1971) derive this same spectral decomposition by way of the general form for the spectral decomposition of a bidiagonal matrix (i.e., a general "pure-birth" process). To apply the spectral decomposition, we let $\mathbf{v}_k$ denote the $k$th row of the eigenvector matrix and we let $\mathbf{w}_k$ denote the $k$th column of the inverse eigenvector matrix. We then have:

$$\begin{aligned}
\mathbb{P}(K_{\acute{n}+n} = k | K_{\acute{n}} = t) = [\mathbf{P}^n]_{t,k} &= (\mathbf{v}\boldsymbol{\Lambda}^n\mathbf{w})_{t,k} \\
&= \mathbf{v}_t\boldsymbol{\Lambda}^n\mathbf{w}_k \\
&= \sum_{i=0}^{m} \lambda_i^n \cdot v_{t,i} \cdot w_{i,k} \\
&= \sum_{i=0}^{m} \left(1 - \theta \cdot \frac{m-i}{m}\right)^n \cdot (-1)^{i-k} \binom{m-t}{i-t} \cdot \binom{m-i}{k-i} \\
&= \sum_{i=0}^{m} \left(1 - \theta \cdot \frac{m-i}{m}\right)^n \cdot (-1)^{i-k} \binom{m}{k} \cdot \binom{k}{i} \cdot \frac{(i)_t}{(m)_t} \\
&= \binom{m}{k} \sum_{i=t}^{k} \binom{k}{i} (-1)^{i-k} \cdot \left(1 - \theta \cdot \frac{m-i}{m}\right)^n \cdot \frac{(i)_t}{(m)_t}.
\end{aligned}$$

(Here we use the "falling factorials" $(m)_t = \prod_{i=0}^{t-1}(m-i)$ to expand binomial coefficients.) This gives a general form for the conditional probabilities in the extended occupancy problem. Taking $t = 0$ gives the marginal form:

$$\begin{aligned}
\mathbb{P}(K_n = k) = [\mathbf{P}^n]_{0,k} &= \binom{m}{k} \sum_{i=0}^{k} \binom{k}{i} (-1)^{i-k} \cdot \left(1 - \theta \cdot \frac{m-i}{m}\right)^n \\
&= \frac{\theta^n}{m^n} \cdot (m)_k \cdot S\left(n, k, m \cdot \frac{1-\theta}{\theta}\right),
\end{aligned}$$

where $S(n, k, \phi)$ are the noncentral Stirling numbers of the second kind (see Appendix I). This gives us a succinct form for the marginal probabilities arising in the occupancy problem.

In the above working, we offered the conditional form of the extended occupancy distribution as our most general form, with the marginal form occurring for $t = 0$. However, it is possible to rewrite the conditional form of the extended occupancy distribution using the marginal form with a corresponding variation in the number of bins and the probability of "falling through" a bin. To see this, we first note that —with a bit of algebra— it can be shown that:

$$\binom{m}{k}\binom{k}{i} \cdot \frac{(i)_t}{(m)_t} = \binom{m-t}{k-t}\binom{k-t}{i}.$$

Hence, an alternative form for the conditional occupancy probability is:

$$\begin{aligned}
\mathbb{P}(K_{\acute{n}+n} = k | K_{\acute{n}} = t) &= \binom{m}{k} \sum_{i=t}^{k} \binom{k}{i} (-1)^{i-k} \cdot \left(1 - \theta \cdot \frac{m-i}{m}\right)^n \cdot \frac{(i)_t}{(m)_t} \\
&= \binom{m-t}{k-t} \sum_{i=t}^{k} \binom{k-t}{i} (-1)^{i-k} \cdot \left(1 - \theta \cdot \frac{m-i}{m}\right)^n \\
&= \binom{m-t}{k-t} \sum_{i=0}^{k-t} \binom{k-t}{t+i} (-1)^{i-k+t} \cdot \left(1 - \theta \cdot \frac{m-i-t}{m}\right)^n \\
&= \binom{m-t}{k-t} \sum_{i=0}^{k-t} \binom{k-t}{t+i} (-1)^{i-k+t} \cdot \left(1 - \theta \cdot \frac{m-t}{m} \cdot \frac{m-t-i}{m-t}\right)^n \\
&= \binom{m-t}{k-t} \sum_{i=0}^{k-t} \binom{k-t}{t+i} (-1)^{i-k+t} \cdot \left(1 - \theta \left(1 - \frac{t}{m}\right) \cdot \frac{m-t-i}{m-t}\right)^n \\
&= \frac{\left(\theta(1-t/m)\right)^n}{(m-t)^n} \cdot (m-t)_{k-t} \cdot S\left(n, k-t, m-t \cdot \frac{1-\theta-t/m}{\theta}\right).
\end{aligned}$$

This is the same form as the marginal occupancy probability given above, but with an argument value of $k-t$ occupied bins out of $m-t$ bins, and with probability parameter $\theta(1-t/m)$.

This result also follows simple intuition. To convert the conditional occupancy probability to a marginal occupancy probability, we can treat the problem as a marginal occupancy problem where we allocate all the new balls to $m-t$ unoccupied bins, but we also have to reduce the probability parameter so that a new ball is considered to "fall through" its bin if it would have been allocated to one of the $t$ bins already occupied by previous balls. In this method, the conditional probability of allocating a new ball to an already occupied bin is "folded into" the probability parameter, allowing us to write the conditional occupancy probability as a marginal occupancy probability. This result shows that both the marginal occupancy probability and the conditional occupancy probability can be written in the same distributional form, which can be used to solve the extended occupancy problem. In the next section, we will formalise this and look at three distributions that arises in our analysis.

Our above analysis gives solutions to the extended occupancy problem. We can make these solutions clearer and more succinct by introducing a class of distributions for the extended occupancy distribution, and naming its parameters.

**Definition (The extended occupancy distribution):** This is a discrete distribution with probability mass function given by:[3]

$$\text{Occ}(k|n,m,\theta) = \frac{\theta^n}{m^n} \cdot (m)_k \cdot S\left(n,k,m \cdot \frac{1-\theta}{\theta}\right) \qquad 0 \leq k \leq \min(n,m)\,,$$

where $m \in \bar{\mathbb{N}}$ is the **space parameter** (number of bins),[4] $n \in \mathbb{N}$ is the **size parameter** (number of balls), and $0 < \theta \leq 1$ is the **probability parameter**.[5] In the special case where $\theta = 1$ the distribution reduces to the classical occupancy distribution $\text{Occ}(k|n,m) \equiv \text{Occ}(k|n,m,1)$. □

Our above analysis gives solutions to the extended occupancy problem, and shows that both the marginal and conditional distributions that arise in this problem are the same distributional form (but with different parameters). Using the notation introduced in our definition of the extended occupancy distribution, the marginal and conditional probabilities of interest are:

$$\mathbb{P}(K_n = k) = \text{Occ}(k|n,m,\theta),$$

$$\mathbb{P}(K_{\acute{n}+n} = k|K_{\acute{n}} = t) = \text{Occ}(k-t|n,m-t,\theta(1-t/m)).$$

The special case where $\theta = 1$ leads to the classical occupancy distribution for the marginal distribution, and in this case the distribution can be derived by a combinatorial argument (see e.g., O'Neill 2021). The main value of extending the classical occupancy distribution to the extended occupancy distribution is that the latter is "closed under conditioning", by which we mean that this family accommodates both the marginal and conditional distributions of the occupancy number. The extended occupancy distribution has been examined by a number of authors including Park (1972), Johnson and Kotz (1977, Section 3.3, pp. pp.139-146), Samuel-Cahn (1974) and Holst (1986). Broader extension to general occupancy problems and corresponding distributions can be found in Charalambides (2005).

---

[3] It is worth noting that we can allow argument values $k > \min(n,m)$ in this expression, and in these cases, the stated expression for the mass function reduces down to zero (see e.g., Ruiz 1996). This means that we can validly use this mass function for all $k = 0,1,2,\ldots$ and the argument values above the specified support will have zero probability. This is a useful aspect of the mass function, since it allows us to use this function for any non-negative argument value, which allows us to play "fast and loose" with the upper bound of sums over the mass function. This will be especially valuable when we look at expected values of functions of the occupancy number (e.g., the moments of the distribution).

[4] We will examine the special case where $m = \infty$ in the next section, and there we will show that the extended occupancy distribution degenerates to the binomial distribution in this case.

[5] In the trivial case where $\theta = 0$ we have $\text{Occ}(k|n,m,\theta) = \mathbb{I}(k=0)$ so that the distribution is a point mass on $k = 0$. This reflects the fact that setting the probability parameter to zero means that all the balls fall through their allocated bins with probability one. This trivial case is not particularly interesting, so we have removed it from the scope of our analysis in this paper. In all subsequent equations we will assume that $\theta > 0$.

**Remark 1:** Mathematically, the mass function of the extended occupancy distribution arises from the expansion $(m+\phi)^n = \sum_{k=0}^{n} (m)_k \cdot S(n,k,\phi)$ for the non-central Stirling numbers of the second kind (see Appendix I). Each of the terms in this sum is non-negative, and those terms constitute a kernel for the mass function of the extended occupancy distribution. □

The occupancy distribution provides a natural extension to the binomial distribution, insofar as it takes the count of the "effective balls" and "squashes" this number down to the occupancy number by counting only those effective balls that are not duplicating the occupancy of a bin. In the next section we will see that the occupancy distribution actually provides a generalisation of the binomial distribution. However, for the moment, it is worth comparing the forms of the mass functions of the two distributions. To do this, it is quite useful to write the mass function of the occupancy distribution in an alternative form as a product of the binomial mass function multiplied by an adjustment term involving the scaled Stirling function (see Appendix I):

$$\begin{aligned}
\text{Occ}(k|n,m,\theta) &= \frac{\theta^n}{m^n} \cdot (m)_k \cdot S\left(n,k,m \cdot \frac{1-\theta}{\theta}\right) \\
&= \frac{(m)_k}{m^k} \cdot \frac{S\left(n,k,m \cdot \frac{1-\theta}{\theta}\right)}{\binom{n}{k}\left(m \cdot \frac{1-\theta}{\theta}\right)^{n-k}} \cdot \binom{n}{k} \cdot \theta^k \cdot (1-\theta)^{n-k} \\
&= \frac{(m)_k}{m^k} \cdot \Pi\left(n,k,m \cdot \frac{1-\theta}{\theta}\right) \times \text{Bin}(k|n,\theta).
\end{aligned}$$

This form shows a close resemblance between the mass functions of these two distributions. Moreover, the monotonicity and limit properties of the scaled Stirling function (see Lemma 1 in Appendix I) allow us to obtain useful properties of the occupancy distribution.

A full account of the properties of the extended occupancy distribution is beyond the scope of this paper. Nevertheless, it is worth giving some basic properties including its moments and asymptotic form, since these are useful for computational purposes. As with many discrete distributions, the moments of the extended occupancy distribution are simplest when presented through the factorial moments. These factorial moments yield corresponding functions for the raw and central moments, which can be computed with a reasonable amount of algebra. We will go as far as the kurtosis of the distribution, noting that the form of this moment is already quite cumbersome. We will also show the asymptotic form of important moments. Higher-order raw and central moments can be computed from the factorial moments, but they are not particularly illuminating.

**THEOREM 3 (Factorial and raw moments):** Letting $E_r \equiv (1 - \theta r/m)^n$, we have:

$$\mathbb{E}((m - K_n)_r) = (m)_r \cdot E_r.$$

We can see from Theorem 3 that the occupancy distribution gives a simple form for the factorial moments, in terms of the terms $E_r$. (This notation comes in handy below when we write the central moments of the distribution.) Since $(m - K_n)_r$ is a polynomial in $K_n$, the factorial moments are used to derive the raw and central moments. The algebra is cumbersome, so for brevity we state them here as corollaries to this above theorem, without further derivation.

**COROLLARY (Central moments):** The extended occupancy distribution has mean, variance, skewness and kurtosis given respectively by:

$$\mu_{n,m,\theta} = m(1 - E_1),$$

$$\sigma^2_{n,m,\theta} = m[(m-1)E_2 + E_1 - mE_1^2],$$

$$\gamma_{n,m,\theta} = -\frac{E_1 - 3E_2 + 2E_3 + m\begin{pmatrix} 3(E_2 - E_1^2) \\ +2m(E_1^3 - E_1 E_2) \\ +(m-3)(E_3 - E_1 E_2) \end{pmatrix}}{m^{1/2}[(m-1)E_2 + E_1 - mE_1^2]^{3/2}},$$

$$\kappa_{n,m,\theta} = \frac{\begin{pmatrix} E_1 - 4mE_1^2 + 6m^2E_1^3 - 3m^3E_1^4 \\ +7(m-1)E_2 + 6(m-1)(m-2)E_3 \\ +(m-1)(m-2)(m-3)E_4 \\ -12m(m-1)E_1E_2 + 6m^2(m-1)E_1^2E_2 \\ -4m(m-1)(m-2)E_1E_3 \end{pmatrix}}{m[(m-1)E_2 + E_1 - mE_1^2]^2}.$$

**COROLLARY (Asymptotic central moments):** As $n \to \infty$ we have the asymptotic equivalence $E_r \sim e^{-\theta r n/m}$ which gives the asymptotic forms:

$$\mu_{n,m,\theta} \sim m(1 - e^{-\theta n/m}),$$

$$\sigma^2_{n,m,\theta} \sim m e^{-\theta n/m}(1 - e^{-\theta n/m}),$$

$$\gamma_{n,m,\theta} \sim -\frac{1}{\sqrt{m}} \cdot \frac{1 - 2e^{-\theta n/m}}{\sqrt{e^{-\theta n/m}(1 - e^{-\theta n/m})}},$$

$$\kappa_{n,m,\theta} \sim 3 + \frac{1}{m} \cdot \frac{1 - 6e^{-\theta n/m}(1 - e^{-\theta n/m})}{e^{-\theta n/m}(1 - e^{-\theta n/m})}.$$

If $n \to \infty$ and $m \to \infty$ in a way that yields a fixed finite limit for $n/m$ then we have $\gamma_{n,m} \to 0$ and $\kappa_{n,m} \to 3$ (so the distribution is asymptotically unskewed and mesokurtic). □

**COROLLARY (Asymptotic central moments):** As $m \to \infty$ we have the asymptotic equivalence $m^a E_r^b \sim \sum_{i=0}^{a} \binom{bn}{i} (-1)^i (r\theta)^i m^{a-i}$ which gives the asymptotic forms:

$$\mu_{n,m,\theta} \sim \theta,$$

$$\sigma_{n,m,\theta}^2 \sim n\theta(1-\theta),$$

$$\gamma_{n,m,\theta} \sim \frac{1}{\sqrt{n}} \cdot \frac{1-2\theta}{\sqrt{\theta(1-\theta)}},$$

$$\kappa_{n,m,\theta} \sim 3 + \frac{1}{n} \cdot \frac{1-6\theta(1-\theta)}{\theta(1-\theta)}.$$

The moments of the extended occupancy distribution give us a reasonable sense of the shape of the distribution. In particular, we see that —under broad limit conditions— the distribution is asymptotically unskewed and mesokurtic. In fact, this is just a partial aspect of a powerful limit result for general occupancy distributions given in Hwang and Janson (2008). If $n \to \infty$ and $m \to \infty$ in such a way that $\sigma_{n,m}^2 \to \infty$ (having a fixed finite limit for $n/m$ is a sufficient condition for this convergence) then the mass function for the extended occupancy distribution converges uniformly to the normal density with the same mean and variance. This result can be used to approximate the occupancy distribution for large values of $n$ and $m$ where it is not feasible to compute the distribution (due to difficulties computing the Stirling numbers of the second kind for large input values).

## 3. The extended occupancy distribution generalises the binomial distribution

Here we will show that the occupancy distribution provides a generalisation of the binomial distribution (i.e., it subsumes the binomial as a special case), which will later allow us to frame various properties and mixture characterisations as extensions of well-known characterisations for the binomial distribution. The attentive reader may already have noticed that our definition of the extended occupancy distribution allows the space parameter (number of bins) to occur in the extended domain $m \in \bar{\mathbb{N}}$, which allows the value $m = \infty$. In this latter case we define the distribution by its limit $\text{Occ}(k|n, \infty, \theta) \equiv \lim_{m\to\infty} \text{Occ}(k|n, m, \theta) = \text{Bin}(k|n, \theta)$. (This case forms part of our definition of the extended occupancy distribution, but we omitted it from the above definition in order to discuss it in more detail here.)

Harkness (1969) notes the similarity of the occupancy distribution to the binomial distribution (pp. 112-114). We will prove the special case identified above formally in a moment, but the simplest way to see that this limit gives the binomial distribution is to observe that the transition probability matrix converges to the infinite dimensional matrix:

$$\mathbf{P}_\infty \equiv \lim_{m\to\infty} \mathbf{P} \equiv \begin{bmatrix} 1-\theta & \theta & 0 & \cdots & 0 & 0 & \cdots \\ 0 & 1-\theta & \theta & \cdots & 0 & 0 & \cdots \\ 0 & 0 & 1-\theta & \cdots & 0 & 0 & \cdots \\ \vdots & \vdots & \vdots & \ddots & \vdots & \vdots & \ddots \\ 0 & 0 & 0 & \cdots & 1-\theta & \theta & \cdots \\ 0 & 0 & 0 & \cdots & 0 & 1-\theta & \cdots \\ \vdots & \vdots & \vdots & \ddots & \vdots & \vdots & \ddots \end{bmatrix}.$$

The reader will recognise this matrix as the transition probability matrix for a Bernoulli process, and its powers give matrices with elements:

$$[\mathbf{P}_\infty^n]_{t,k} \equiv \text{Bin}(k-t|n,\theta).$$

Thus, starting in state $K_0 = 0$ and taking $n$ steps (i.e., allocating $n$ balls at random) gives the state $K_n$ having a binomial distribution with size parameter $n$ and probability parameter $\theta$. We have derived the result heuristically here, but it can be formally established either via analysis of the limiting properties of Markov chains, or by purely algebraic analysis on the established probability mass function for the extended occupancy distribution.

**THEOREM 4A (Generalisation of the binomial distribution):** The occupancy distribution satisfies the limiting form:

$$\lim_{m\to\infty} \text{Occ}(k|n,m,\theta) = \text{Bin}(k|n,\theta).$$

This theorem formally establishes that the extended occupancy distribution is a generalisation of the binomial distribution. Intuitively, the limiting result reflects the fact that, with infinite bins, there is zero probability that any two balls will fall in the same bin. Thus, the occupancy number is then the "effective" number of balls that have not "fallen through" their allocated bins, which is merely a count of independent Bernoulli random variables with fixed probability. Indeed, going back to our initial setup for the occupancy problem, we note that the number of effective balls in the problem is $n_{\text{eff}} = \sum_{i=1}^n \mathbb{I}(U_i \neq \bullet) = \sum_{i=1}^n \mathbb{I}(Q_i = 1)$, with the underlying values $Q_1, Q_2, Q_3, \ldots \sim \text{IID Bern}(\theta)$.

Since the occupancy distribution generalises the binomial, and since the binomial distribution is well-known to lead to other distributional forms (e.g., the Poisson) under appropriate limits, a natural follow-up question is to ask whether we obtain any interesting distribution if we keep the space parameter $m$ as a finite value, but take the limits on the other parameters that would yield the Poisson distribution from the binomial. It turns out that this limiting exercise leads us to another binomial distribution (with different parameters). Since the binomial is itself a generalisation of the Poisson distribution, under appropriate limits, the occupancy distribution generalises both the binomial and the Poisson.[6]

**THEOREM 4B (Limit to the binomial distribution):** The occupancy distribution satisfies the limiting form:

$$\lim_{\substack{n\to\infty,\theta\to 0\\ n\theta\to\lambda}} \text{Occ}(k|n,m,\theta) = \text{Bin}\left(k\middle|m, 1-\exp\left(-\frac{\lambda}{m}\right)\right).$$

The binomial distribution has a number of well-known properties relating to recursion on its parameters , monotone likelihood ratio with respect to its parameters, and resulting stochastic dominance (see e.g., Johnson and Kotz 1969, pp. 50-86). Since the occupancy distribution generalises the binomial, it is useful to see how the properties of the binomial are generalised in the occupancy distribution. In particular, it is well-known that a binomial random variable has a monotone likelihood ratio and is therefore "stochastically increasing" in $n$ and $\theta$ (in the sense of first-order stochastic dominance), and that it obeys simple recursive and differential equations on its parameters. We now establish a set of generalised equations, monotonicity and stochastic dominance results that apply to the extended occupancy distribution, with the binomial recursive equations and stochastic dominance results occurring as a special case.

**THEOREM 5 (Recursive and differential equations):** The extended occupancy distribution satisfies the following recursive/differential equations:

$$\text{Occ}(k|n+1,m,\theta) = \theta\cdot\frac{m-k+1}{m}\cdot\text{Occ}(k-1|n,m,\theta) + \left(1-\theta\cdot\frac{m-k}{m}\right)\cdot\text{Occ}(k|n,m,\theta),$$

$$\text{Occ}(k|n,m+1,\theta) = \frac{m+1}{m-k+1}\cdot\left(1-\frac{\theta}{m+1}\right)^n\cdot\text{Occ}\left(k\middle|n,m,\frac{m\theta}{1-\theta+m}\right),$$

$$\frac{\partial}{\partial\theta}\text{Occ}(k|n,m,\theta) = -\frac{m-k}{m}\cdot\text{Occ}(k|n-1,m,\theta) + \frac{m-k+1}{m}\cdot\text{Occ}(k-1|n-1,m,\theta).$$

[6] The latter is easily achieved by taking $n\to\infty$, $m\to\infty$ and $\theta\to 0$ such that $n\theta\to\lambda$. Taking these limits gives $m(1-\exp(-\lambda/m))\to\lambda$ so that $\text{Occ}(k|n,m,\theta)\to\text{Pois}(k|\lambda)$.

**COROLLARY:** In the case where $m = \infty$ we have the binomial recursive/differential equations:

$$\text{Bin}(k|n+1,\theta) = \theta \cdot \text{Bin}(k-1|n,\theta) + (1-\theta) \cdot \text{Bin}(k|n,\theta),$$

$$\frac{\partial}{\partial\theta}\text{Bin}(k|n,\theta) = -\text{Bin}(k|n-1,\theta) + \text{Bin}(k-1|n-1,\theta).$$

**THEOREM 6 (First-order stochastic dominance):** Let $F(k|n,m,\theta) \equiv \mathbb{P}(K_n \leq k)$ denote the cumulative distribution function for the extended occupancy distribution. This satisfies the following first-order stochastic dominance relations:

$$n \leq n' \quad \Longrightarrow \quad F(k|n,m,\theta) \geq F(k|n',m,\theta) \quad \text{(strict if } n < n' \text{ and } m > 1\text{)},$$

$$m \leq m' \quad \Longrightarrow \quad F(k|n,m,\theta) \geq F(k|n,m',\theta) \quad \text{(strict if } m < m' \text{ and } n > 1\text{)},$$

$$\theta \leq \theta' \quad \Longrightarrow \quad F(k|n,m,\theta) \geq F(k|n,m,\theta') \quad \text{(strict if } \theta < \theta'\text{)}.$$

Theorem 5 establishes equations for the occupancy distribution generalising similar equations for the binomial distribution. Theorem 6 shows that the stochastic dominance results for the binomial distribution also hold for the extended occupancy distribution, and are also extended to the new parameter $m$. The theorem establishes that an extended occupancy random variable is "stochastically increasing" in $n$, $m$ and $\theta$ (in the sense of first-order stochastic dominance). This result accords with intuition, since increasing the number of balls, the number of bins, or the probability of occupancy, will all tend to increase the number of occupied bins. This stochastic dominance result also gives us a useful intuitive understanding of the effect of the generalisation from the binomial to the occupancy distribution. By imposing a finite parameter $m$ (rather than the value $m = \infty$ that gives the binomial distribution) we "squash" the effective balls into a finite number of available bins, which gives rise to the possibility that more than one effective ball will share a bin, so the excess balls will not count towards our occupancy number. A simple corollary of Theorem 6 is that imposition of a finite value of $m$ (instead of $m = \infty$) will tend to reduce the value of the occupancy number.

In Appendix II we prove the stochastic dominance results in Theorem 6 algebraically from the mass function for the occupancy distribution. However, a simpler intuition for the result is obtained by noting that the Markov chain defining the extended occupancy distribution has a number of monotonicity properties. The rows of the matrix have cumulative sums that are increasing, so the matrix is "monotone" in the sense described in Daley (1968).[7] Moreover, in any row of the matrix, the cumulative sum of terms is decreasing in $m$ and $\theta$ (in the degenerate case where $\theta = 0$ it is only non-increasing in $m$). An alternative proof (not pursued here) could be couched in terms of these general properties of monotone Markov chains.

## 4. Excess hitting times and the negative occupancy distribution

Another distribution related to the binomial distribution is the negative binomial distribution, and it turns out that we can find a natural extension of this latter distribution arising in the occupancy process. This distribution is obtained by considering the "excess hitting time"[8] for the event $K_n = k$ in the Markov chain, which we denote by:

$$T_k \equiv \min\{t = 0,1,2,\dots | K_{k+t} = k\}.$$

Derivation of the distribution of this excess hitting time is quite straight-forward. The event $T_k \leq t$ is equivalent to the event $K_{k+t} \geq k$. Hence, for all $0 < k \leq m$ and $t \geq 0$ the cumulative distribution of $T_k$ can be obtained from the occupancy distribution as:

$$F_{T_k}(t) = \mathbb{P}(K_{k+t} \geq k) = \sum_{r=0}^{m-k} \text{Occ}(k + r | k + t, m, \theta).$$

To find the probability mass function we take the first difference of the distribution function. As a preliminary step, applying the first recursive equation in Theorem 4 we obtain:

$$\begin{aligned} \text{Occ}(k + r | k + t, m, \theta) &= \theta \cdot \frac{m - k - r + 1}{m} \cdot \text{Occ}(k + r - 1 | k + t - 1, m, \theta) \\ &\quad + \left(1 - \theta \cdot \frac{m - k - r}{m}\right) \cdot \text{Occ}(k + r | k + t - 1, m, \theta). \end{aligned}$$

We therefore have:

$$\begin{aligned} &\text{Occ}(k + r | k + t, m, \theta) - \text{Occ}(k + r | k + t - 1, m, \theta) \\ &\qquad = \frac{\theta}{m} \begin{bmatrix} (m - k - r + 1) \cdot \text{Occ}(k + r - 1 | k + t - 1, m, \theta) \\ -(m - k - r) \cdot \text{Occ}(k + r | k + t - 1, m, \theta) \end{bmatrix}. \end{aligned}$$

[7] That is, the function $\Lambda(t, k) \equiv \mathbb{P}(K_{n+1} \leq k | K_n = t)$ is strictly increasing in $t$ for every $k$. (In the degenerate case where $\theta = 0$ it is only non-decreasing.) This can be established from the transition probability matrix **P**.

[8] Since the occupancy process is a "pure birth process" the occupancy number can increase by no more than one unit with each ball, which means that the occupancy number $K_n = k$ requires at least $n = k$ balls. Our "excess" hitting time is measuring the number of "excess" balls —i.e., the number of balls beyond this minimum— to obtain the stipulated occupancy number.

Using this result, the mass function for the excess hitting time is:

$$\begin{aligned}
\mathbb{P}(T_k = t) &= F_{T_k}(t) - F_{T_k}(t-1) \\
&= \sum_{r=0}^{m-k} [\text{Occ}(k+r|k+t,m,\theta) - \text{Occ}(k+r|k+t-1,m,\theta)] \\
&= \frac{\theta}{m} \cdot \left[ \begin{array}{c} \sum_{r=0}^{m-k} (m-k-r+1) \cdot \text{Occ}(k+r-1|k+t-1,m,\theta) \\ - \sum_{r=1}^{m-k+1} (m-k-r+1) \cdot \text{Occ}(k+r-1|k+t-1,m,\theta) \end{array} \right] \\
&= \frac{\theta}{m} \cdot \left[ \begin{array}{c} (m-k-r+1) \cdot \text{Occ}(k+r-1|k+t-1,m,\theta)|_{r=0} \\ -(m-k-r+1) \cdot \text{Occ}(k+r-1|k+t-1,m,\theta)|_{r=m-k+1} \end{array} \right] \\
&= \frac{\theta}{m} \cdot (m-k+1) \cdot \text{Occ}(k-1|k+t-1,m,\theta) \\
&= \theta \cdot \frac{m-k+1}{m} \cdot \frac{\theta^{k+t-1}}{m^{k+t-1}} \cdot (m)_{k-1} \cdot S\left(k+t-1, k-1, m \cdot \frac{1-\theta}{\theta}\right) \\
&= \frac{\theta^{k+t}}{m^{k+t}} \cdot (m)_k \cdot S\left(k+t-1, k-1, m \cdot \frac{1-\theta}{\theta}\right).
\end{aligned}$$

Below we formalise this mass function as defining a class of distributions we call the "negative occupancy distribution". We also consider a special case of this distribution, which describes the behaviour of the excess coupons collected in the famous coupon-collector problem.

**DEFINITION (The negative occupancy distribution):** This is a discrete distribution that has probability mass function given over all integer arguments $t \geq 0$ as:

$$\text{NegOcc}(t|m,k,\theta) = \frac{\theta^{k+t}}{m^{k+t}} \cdot (m)_k \cdot S\left(k+t-1, k-1, m \cdot \frac{1-\theta}{\theta}\right),$$

where $0 < k \leq m \leq \infty$ are the **occupancy parameter** (number of occupied bins) and **space parameter** (number of bins) respectively, and $0 < \theta \leq 1$ is the **probability parameter**. □

**DEFINITION (Coupon-collector distribution):** This is a discrete distribution with probability mass function given over all integer argument values $t \geq 0$ as:

$$\text{CoupColl}(t|m,\theta) = \frac{\theta^{m+t}}{m^{m+t}} \cdot m! \cdot S\left(m+t-1, m-1, m \cdot \frac{1-\theta}{\theta}\right).$$

where $0 < m < \infty$ is the **space parameter** (number of bins) and $0 < \theta \leq 1$ is the **probability parameter**. Since $\text{CoupColl}(t|m,\theta) = \text{NegOcc}(t|m,m,\theta)$ the coupon-collector distribution is a special case of the negative occupancy distribution with $k = m$. □

**REMARK 2:** Mathematically, the mass function of the negative occupancy distribution arises from the ordinary generating function for the noncentral Stirling numbers of the second kind, which is $\sum_{n=k}^{\infty} S(n,k,\phi) \cdot x^n = x^k / \prod_{i=0}^{k} (1 - (i + \phi)x)$ (see Appendix I). Substituting the noncentrality parameter $\phi = m\,(1-\theta)/\theta$ and the argument $x = \theta/m$ and rearranging gives the norming equation for the mass function of the negative occupancy distribution. □

Since $\mathbb{P}(T_k = t) = \text{NegOcc}(t|m,k,\theta)$ we can see that the negative occupancy distribution is the appropriate family of distributions to describe the behaviour of the excess hitting time in the occupancy process. In the case where $k = m$ we are looking at the excess number of balls that are required to fully occupy all the bins in the occupancy problem. In this case, we have referred to the distribution as the "coupon-collector" distribution. The classical version of the coupon-collector distribution arises in the "coupon-collector problem", which examines the number of randomly obtained coupons that need to be collected to obtain a full set (see e.g., Dawkins 1991, Adler et al 2003). Note that our distribution describes the *excess* number of coupons required for a full set, not the *total* number of required coupons; it is trivial to convert to the distribution of the total required number of coupons if required.

As with the occupancy distribution above, it is useful to write the mass function of the negative occupancy distribution as a product of the negative binomial mass function multiplied by an adjustment term involving the scaled Stirling function (Appendix I):

$$\begin{aligned}
\text{NegOcc}(t|m,k,\theta) &= \frac{\theta^{k+t}}{m^{k+t}} \cdot (m)_k \cdot S\left(k+t-1, k-1, m \cdot \frac{1-\theta}{\theta}\right) \\
&= \frac{(m)_k}{m^k} \cdot \frac{S\left(k+t-1, k-1, m \cdot \frac{1-\theta}{\theta}\right)}{\binom{k+t-1}{k-1}\left(m \cdot \frac{1-\theta}{\theta}\right)^t} \cdot \binom{k+t-1}{t} \cdot (1-\theta)^t \cdot \theta^k \\
&= \frac{(m)_k}{m^k} \cdot \Pi\left(k+t-1, k-1, m \cdot \frac{1-\theta}{\theta}\right) \times \text{NegBin}(t|k, 1-\theta).
\end{aligned}$$

This alternative form shows that there is a close resemblance between the mass functions of these two distributions.[9] In fact, it is simple to show that the negative occupancy distribution

[9] In this equation we are using the "standard" parameterisation of the negative binomial distribution, giving the probability of seeing $t$ "failures" until we first get $k$ "successes", where the probability parameter in the mass function is the probability of a "failure". In case of any doubt, we stipulate the definition:

$$\text{NegBin}(t|k,p) = \binom{k+t-1}{t} \cdot p^t \cdot (1-p)^k \qquad \text{for all } t = 0,1,2,\ldots.$$

generalises the negative binomial distribution. In the case $m = \infty$ we define the distribution by the limit $\text{NegOcc}(t|\infty, k, \theta) \equiv \lim_{m\to\infty} \text{NegOcc}(t|m, k, \theta) = \text{NegBin}(t|k, 1-\theta)$.

**THEOREM 7 (Generalisation of negative binomial distribution):** The negative occupancy distribution satisfies the limiting form:

$$\lim_{m\to\infty} \text{NegOcc}(t|m, k, \theta) = \text{NegBin}(t|k, 1-\theta).$$

The negative occupancy distribution generalises the negative binomial distribution by adding a space parameter $m$ that allows the occupancy to be "squashed" into a finite number of bins.[10] The distribution obeys a number of recursive/differential results, and stochastic dominance properties that generalise results for the negative binomial distribution.

**THEOREM 8 (Recursive and differential equations):** The negative occupancy distribution has recursive/differential equations with respect to its parameters given by:

$$\text{NegOcc}(t|m+1, k, \theta) = \frac{m+1}{m-k+1} \cdot \left(1 - \frac{\theta}{m+1}\right)^{k+t} \cdot \text{NegOcc}\left(t\middle|m, k, \frac{m\theta}{1-\theta+m}\right),$$

$$\text{NegOcc}(t|m, k+1, \theta) = \theta \cdot \frac{m-k}{m} \cdot \sum_{i=0}^{t} \left(1 - \theta \cdot \frac{m-k}{m}\right)^{i} \cdot \text{NegOcc}(t-i|m, k, \theta),$$

$$\frac{\partial}{\partial\theta} \text{NegOcc}(t|m, k, \theta) = \frac{1}{\theta} \cdot \text{NegOcc}(t|m, k, \theta) + \frac{m-k+1}{m} \cdot \begin{bmatrix} \text{NegOcc}(t|m, k-1, \theta) \\ -\text{NegOcc}(t-1|m, k, \theta) \end{bmatrix}.$$

**COROLLARY:** In the case where $m = \infty$ we have the recursive/differential equations for the negative binomial distribution:

$$\text{NegBin}(t|k+1, 1-\theta) = \theta \cdot \sum_{i=0}^{t} (1-\theta)^i \cdot \text{NegBin}(t-i|k, 1-\theta),$$

$$\frac{\partial}{\partial\theta} \text{NegBin}(t|k, 1-\theta) = \frac{1}{\theta} \cdot \text{NegBin}(t|k, 1-\theta) + \begin{bmatrix} \text{NegBin}(t|k-1, 1-\theta) \\ -\text{NegBin}(t-1|k, 1-\theta) \end{bmatrix}.$$

[10] Note that the negative binomial distribution is parameterised in terms of the probability $1-\theta$ of falling through the bin. The negative binomial counts the number of "failures" before a given number of "successes", where the distribution is parameterised in terms of the probability of "failure". In the context of our analysis the "failures" balls that do not contribute to the occupancy number (either because they fall through their bins or because they occupy a bin that is already occupied) and the "successes" are balls that occupy a new bin and therefore increase the occupancy number. In the case where $m = \infty$ the parameter $\theta$ is the probability of occupancy, which is the probability of a "success". The corresponding probability of a "failure" is $1-\theta$, which is why this value appears as the probability parameter for the negative binomial distribution.

Theorem 8 gives recursive/differential equations for the negative occupancy distribution. The corollary shows that these equations are extensions of well-known equations for the negative binomial distribution. (The reader should note that corresponding recursive equations for the coupon-collector distribution are a little more complicated than for the negative occupancy distribution, owing to the fact that two parameters are collapsed into one. Detailed analysis of this distribution is outside the scope of the present paper.) These recursive equations give rise to corresponding stochastic dominance results, as shown in the theorem below.

**THEOREM 9 (First-order stochastic dominance):** Let $F(t|m,k,\theta) \equiv \mathbb{P}(T_k \leq t)$ denote the cumulative distribution function for the negative occupancy distribution. This satisfies the following first-order stochastic dominance relations:

$$
\begin{array}{llll}
m \leq m' & \Longrightarrow & F(t|m,k,\theta) \leq F(t|m',k,\theta) & \begin{pmatrix} \text{strict if } m < m' \text{ and} \\ k > 1 \text{ or } \theta < 1 \end{pmatrix} \\
k \leq k' & \Longrightarrow & F(t|m,k,\theta) \geq F(t|m,k',\theta) & (\text{strict if } k < k') \\
\theta \leq \theta' & \Longrightarrow & F(t|m,k,\theta) \leq F(t|m,k,\theta') & (\text{strict if } \theta < \theta')
\end{array}
$$

In the proofs in Appendix II, we establish the above stochastic dominance results algebraically from the mass function of the negative occupancy distribution. These results follow directly from the monotone likelihood-ratio properties of the distribution with respect to its parameters, but they also have some basic statistical intuition. In particular, the stochastic dominance results for the negative occupancy distribution are intuitively related to stochastic dominance for the extended occupancy distribution —*ceteris paribus*, increases in either $m$ or $\theta$ will tend to increase the occupancy number for any fixed number of balls, and will thus tend to decrease the number of excess balls required to achieve a fixed occupancy number. Contrarily, if we increase the occupancy number $k$ this will tend to directly increase the excess hitting time, since the hitting time is now for a larger outcome value in a pure birth process.

The negative occupancy distribution provides us with a description of the stochastic behaviour of the excess number of balls required to achieve a given occupancy number in the extended occupancy problem. By a simple shift in location, it can also be used to describe the stochastic behaviour of the total number of balls required to achieve a given occupancy number. The coupon-collector distribution is a special case of the negative occupancy distribution, which solves the famous "coupon-collector problem", giving a full description of the behaviour of the minimum number of balls needed to achieve full occupancy.

## 5. The "spillage" and its conditional distribution

From our previous analysis, we have already seen that taking $m = \infty$ means that each ball falls into a different bin, yielding standard Bernoulli sampling. In this special case the occupancy number must be equal to the effective number of balls in the occupancy problem, and thus, we will have $n_{\text{eff}} - K_n = 0$. If we use a finite number of bins $m < \infty$ this is no longer guaranteed, since it is possible for some of the effective balls to occupy the same bin, so that the effective number of balls may exceed the occupancy number. If we consider a bin containing a single ball to be occupied, we can consider the value $n_{\text{eff}} - K_n$ to constitute "spillage" of balls in excess of the number required to occupy the occupied bins.

The third distribution we will examine in the extended occupancy problem is the conditional distribution of the "spillage", conditional on the occupancy number. We will derive this distribution using Bayes' theorem. Conditional on the effective number of balls $n_{\text{eff}} = s$, the distribution of the occupancy number is the classical occupancy distribution:

$$\mathbb{P}(K_n = k | n_{\text{eff}} = s, n, m, \theta) = \text{Occ}(k|s, m).$$

Thus, for all occupancy values $1 \leq k \leq \min(n, m)$ and all argument values $k \leq s \leq n$ for the effective number of balls, we have:

$$\begin{aligned}
\mathbb{P}(n_{\text{eff}} = s | K_n = k, n, m, \theta) &= \frac{\mathbb{P}(K_n = k | n_{\text{eff}} = s, n, m, \theta) \times \mathbb{P}(n_{\text{eff}} = s | n, m, \theta)}{\mathbb{P}(K_n = k | n, m, \theta)} \\
&= \frac{\text{Occ}(k|s, m) \times \text{Bin}(s|n, \theta)}{\text{Occ}(k|n, m, \theta)} \\
&= \frac{\frac{1}{m^s} \cdot (m)_k \cdot S(s, k) \times \binom{n}{s} \cdot \theta^s \cdot (1-\theta)^{n-s}}{\frac{\theta^n}{m^n} \cdot (m)_k \cdot S\left(n, k, m \cdot \frac{1-\theta}{\theta}\right)} \\
&= \frac{\frac{\theta^n}{m^n} \cdot (m)_k \cdot S(s, k) \times \binom{n}{s} \cdot \left(m \cdot \frac{1-\theta}{\theta}\right)^{n-s}}{\frac{\theta^n}{m^n} \cdot (m)_k \cdot S\left(n, k, m \cdot \frac{1-\theta}{\theta}\right)} \\
&= \binom{n}{s} \cdot \left(m \cdot \frac{1-\theta}{\theta}\right)^{n-s} \cdot \frac{S(s, k)}{S(n, k, m \cdot (1-\theta)/\theta)}.
\end{aligned}$$

Taking $s = k + r$ gives the conditional distribution of the "spillage", which is:

$$\begin{aligned}
\mathbb{P}(n_{\text{eff}} - K_n = r | K_n = k, n, m, \theta) &= \mathbb{P}(n_{\text{eff}} = k + r | K_n = k, n, m, \theta) \\
&= \binom{n}{k+r} \cdot \left(m \cdot \frac{1-\theta}{\theta}\right)^{n-k-r} \cdot \frac{S(k+r, k)}{S(n, k, m \cdot (1-\theta)/\theta)}.
\end{aligned}$$

**DEFINITION (The spillage distribution):** This distribution is a discrete probability distribution with probability mass function given by:[11]

$$\text{Spillage}(r|n,k,\phi) \equiv \binom{n}{k+r} \cdot \phi^{n-k-r} \cdot \frac{S(k+r,k)}{S(n,k,\phi)} \qquad r = 0,\dots,n-k,$$

where $n \in \mathbb{N}$ is the **size parameter** (number of balls), $0 \leq k \leq n$ is the **occupancy parameter** (occupancy number) and $0 \leq \phi \leq \infty$ is the **scale parameter**. □

We can see that the spillage distribution describes the behaviour of the "spillage" given our knowledge of the occupancy number. The distribution can also be shifted to describe the behaviour of the number of effective balls given our knowledge of the occupancy number. These two conditional probabilities are given respectively by:[12]

$$\mathbb{P}(n_{\text{eff}} - K_n = r | K_n = k, n, m, \theta) = \text{Spillage}\left(r \middle| n, k, m \cdot \frac{1-\theta}{\theta}\right),$$

$$\mathbb{P}(n_{\text{eff}} = s | K_n = k, n, m, \theta) = \text{Spillage}\left(k + s \middle| n, k, m \cdot \frac{1-\theta}{\theta}\right).$$

It is worth noting here that the distribution of the "spillage" depends on $m$ and $\theta$ only through the scale parameter $\phi = m \cdot (1-\theta)/\theta$. In the classical case where $\theta = 1$ we have the scale parameter $\phi = 0$ so $n_{\text{eff}} = n$ with probability one (and the corresponding "spillage" is $n - k$). In the case where $m = \infty$ and $0 < \theta < 1$ we have the scale parameter $\phi = \infty$ so $n_{\text{eff}} = k$ with probability one (and the corresponding "spillage" is zero). We have been unable to identify this distribution in the existing mathematical or statistical literature, and so to our knowledge it is a "new" distributional family; the name we have ascribed here is our own creation. The name we have chosen reflects the fact that the distribution arises when we consider excess balls above what is required to occupy a bin to "spill" over the capacity of the bin.

---

[11] The case where $\phi = 0$ gives the point-mass distribution $\text{Spillage}(r|n,k,0) = \mathbb{I}(r = n-k)$. Additionally, the case where $\phi = \infty$ is defined via the appropriate limit, as the point-mass distribution:

$$\text{Spillage}(r|n,k,\infty) \equiv \lim_{\phi \to \infty} \text{Spillage}(r|n,k,\phi) = \mathbb{I}(r = 0).$$

[12] To avoid ambiguity, we need to comment on the case where $m = \infty$ and $\theta = 1$, which gives an indeterminate form for the scale parameter. This case corresponds to allocation of balls to an infinite number of bins, with zero probability of falling through the bins, and so it should give $n_{\text{eff}} = n = k$ with probability one. By convention, in this case we set $\phi \equiv 0$ so that the distribution for the spillage is a point mass on the value $r = n - k$. Note that even with this convention, the fact that the values $n$ and $k$ are conditioning parameters in the distribution, allows the user to stipulate values $n \neq k$, and this can give pathological results. This occurs because we are dealing with a conditional distribution where it is possible to stipulate conditioning values that occur jointly with probability zero; the fact that an unusual result can occur is of no consequence.

**REMARK 3:** Mathematically, the mass function of the spillage distribution arises from the well-known expansion for the noncentral Stirling numbers of the second kind in terms of the central Stirling numbers of the second kind (see Appendix I), which can be written as:

$$S(n, k, \phi) = \sum_{r=0}^{n-k} \binom{n}{k+r} \cdot \phi^{n-k-r} \cdot S(k+r, k) .$$

For $\phi \geq 0$ the terms in this sum are non-negative and the terms give the kernel of the mass function of the spillage distribution. □

Unlike our previous two distributions, the present distribution does not generalise any common non-trivial distribution arising as a variant of the binomial. In fact, we see from the theorem below that in the case where we have an infinite number of bins (giving $\phi = \infty$) the distribution degenerates down to a point mass on $r = 0$, reflecting the fact that there is no "spillage" in this case. Thus, rather than providing a useful generalisation of an existing distribution, like our previous distributions, the spillage distribution is a new form that describes the divergence between the effective number of balls and the occupancy number in the setting of the extended occupancy problem. In the case of a finite number of bins, the occupancy number may be "squashed" down below the effective number of balls by the fact that balls may share bins with non-zero probability.

**THEOREM 10 (Limit of the spillage distribution with infinite bins):** The spillage distribution satisfies the limiting form:

$$\lim_{\phi \to \infty} \text{Spillage}(r|n, k, \phi) = \mathbb{I}(r = 0).$$

Taking $n = k + r$ gives the probability that the effective number of balls is equal to the full number of balls (i.e., that all balls were effective). The probability of this outcome, written in terms of the occupancy number $k$ and the spillage $r$ is:

$$\text{Spillage}(r|k + r, k, \phi) = \frac{S(k + r, k)}{S(k + r, k, \phi)} .$$

As we did above with our first two distributions, it is useful to write the mass function of the spillage distribution in an alternative form involving the scaled Stirling function (see Appendix I). With a bit of algebra it can be shown that:

$$\binom{n}{k+r} = \frac{\binom{n}{k}\binom{n-k}{r}}{\binom{k+r}{k}}.$$

This gives an alternative form for the mass function of the spillage distribution:

$$\begin{aligned}
\text{Spillage}(r|n,k,\phi) &= \binom{n}{k+r} \cdot \phi^{n-k-r} \cdot \frac{S(k+r,k)}{S(n,k,\phi)} \\
&= \frac{\binom{n}{k}\binom{n-k}{r}}{\binom{k+r}{k}} \cdot \phi^{n-k-r} \cdot \frac{S(k+r,k)}{S(n,k,\phi)} \\
&= \binom{n-k}{r} \cdot \frac{S(k+r,k,\phi)}{\binom{k+r}{k} \cdot \phi^{r}} \Bigg/ \frac{S(n,k,\phi)}{\binom{n}{k} \cdot \phi^{n-k}} \times \frac{S(k+r,k)}{S(k+r,k,\phi)} \\
&= \binom{n-k}{r} \cdot \frac{\Pi(k+r,k,\phi)}{\Pi(n,k,\phi)} \times \text{Spillage}(r|k+r,k,\phi).
\end{aligned}$$

This form of the mass function frames the probabilities relative to the conditional probability that all the balls in the occupancy problem are effective. As can be seen, the form involves writing the mass function as the product of this probability and an adjustment term involving the scaled Stirling numbers.

As with the other two distributions we have examined in this paper, it is possible to use the recursive/differential properties of the noncentral Stirling numbers of the second kind to obtain corresponding recursive/differential equations for the spillage distribution. These equations are rather cumbersome, and not particularly illuminating, so they are omitted here. As should be unsurprising, the spillage is stochastically increasing in $n$ and decreasing in $k$. It is also possible to establish that the spillage is stochastically decreasing in $\phi$, which means it is stochastically decreasing in $m$ and stochastically increasing in $\theta$.

## 6. Mixture properties involving the occupancy distributions

The three occupancy distributions we have examined have a number of interesting mixture characterisations that are useful for computational and analytic purposes. We will look at each of the distributions in the order presented in our previous examination, and derive mixture characterisations for each, beginning with the extended occupancy distribution. The mixtures in Theorems 12-13 below are also shown in Harkness (1969) (Eqs 24 and 23 respectively). To the knowledge of the present author, the remaining results are new.

One way to derive mixture results for the extended occupancy distribution is to treat the number of balls in the occupancy problem as a random variable with a specified distribution over the natural numbers. This leads to a general mixture form shown in the theorem below, where the marginal mass function of the occupancy number involves the probability generating function of the underlying distribution.

**THEOREM 11 (Random number of balls):** Suppose we let the number of balls in the occupancy problem be a random variable $N \sim p_N$ and let $G_N$ be the corresponding probability generating function of $N$. Then the marginal mass function for the occupancy number is:

$$\mathbb{P}(K_N = k|m,\theta) = \binom{m}{k} \sum_{i=0}^{k} \binom{k}{i} (-1)^{k-i} G_N \left(1 - \theta \cdot \frac{m-i}{m}\right).$$

There are some distributions with simple probability generating functions that conform nicely with this sum expression above, in such a way as to yield useful mixture characterisations. In Theorems 12-13 below we give our first mixture results, which show that binomial and Poisson mixtures of the occupancy distribution both give rise to simple marginal distributions. Later we show some mixture results for the negative occupancy and spillage distributions.

**THEOREM 12 (Occupancy distribution is a binomial mixture of occupancy distributions):** The occupancy distribution satisfies the equation:

$$\text{Occ}(k|n,m,\gamma\theta) = \sum_{r=0}^{n} \text{Bin}(r|n,\theta) \cdot \text{Occ}(k|r,m,\gamma).$$

In the special case where $\gamma = 1$ we obtain the useful mixture equation:

$$\text{Occ}(k|n,m,\theta) = \sum_{r=0}^{n} \text{Bin}(r|n,\theta) \cdot \text{Occ}(k|r,m).$$

Theorem 12 has a simple intuition when we interpret it as involving two independent events resulting in "falling through" the bins. The binomial distribution in the mixture gives the number of "effective" balls that do not fall through the bins due to the new event, and each of these terms is multiplied by the occupancy distribution without the probability of that event incorporated. We have stated the theorem in a general form, but the most important case occurs when $\gamma = 1$, which allows us to write the occupancy distribution as a binomial mixture of the classical occupancy distribution. (This mixture is useful for computational purposes; it can be combined with the algorithms in O'Neill 2021 to yield an algorithm to compute the extended occupancy distribution.)

**Theorem 13 (Binomial distribution is a Poisson mixture of occupancy distributions):** The binomial distribution satisfies the equation:

$$\text{Bin}\left(k \middle| m, 1 - \exp\left(-\frac{\lambda\theta}{m}\right)\right) = \sum_{r=0}^{\infty} \text{Pois}(r|\lambda) \cdot \text{Occ}(k|r, m, \theta) .$$

With $\theta > 0$ this can be written to yield a binomial distribution with parameter $0 < \gamma < 1$ as:

$$\text{Bin}(k|m, \gamma) = \sum_{r=0}^{\infty} \text{Pois}\left(r \middle| m \cdot \frac{|\ln(1-\gamma)|}{\theta}\right) \cdot \text{Occ}(k|r, m, \theta) .$$

Theorem 13 gives a mixture characterisation of the binomial distribution as a Poisson mixture of underlying occupancy distributions. We have already noted that the occupancy distribution is a generalisation of the binomial distribution, so this gives us yet another characterisation of the binomial distribution. Theorems 12-13 are extensions of well-known characterisations of the binomial and Poisson distributions. Taking $m \to \infty$ gives $\text{Occ}(k|r, m, \gamma) \to \text{Bin}(k|r, \gamma)$ so that Theorem 12 reduces down to the well-known mixture:

$$\text{Bin}(k|n, \gamma\theta) = \sum_{r=0}^{n} \text{Bin}(r|n, \theta) \cdot \text{Bin}(k|r, \gamma) ,$$

Taking $m \to \infty$ and $\gamma \to 0$ with $m\gamma \to \lambda\theta$ we can apply L'Hôpital's rule to show that:

$$m \cdot \frac{|\ln(1-\gamma)|}{\theta} = \frac{m\gamma}{\theta} \cdot \frac{|\ln(1-\gamma)|}{\gamma} \to \frac{m\gamma}{\theta} \cdot \frac{1}{1-\gamma} \to \lambda.$$

Since these limits give $\text{Occ}(k|r, m, \theta) \to \text{Bin}(k|r, \theta)$ and $\text{Bin}(k|m, \gamma) \to \text{Pois}(k|\lambda\theta)$ we see that the mixture equation in Theorem 13 (second equation) reduces to the well-known mixture:

$$\text{Pois}(k|\lambda\theta) = \sum_{r=0}^{\infty} \text{Pois}(r|\lambda) \cdot \text{Bin}(k|r, \theta) .$$

The extended occupancy distribution provides a useful extension to the binomial distribution, with mixture properties that connect it to other common discrete distributions that arise in statistical practice. In particular, we have established that a Poisson mixture of occupancy distributions yields the binomial distribution, and thus provides a natural link between these two distributions. By viewing the occupancy distribution as a distribution relating to a Markov chain we are able to establish its properties either via the theory of Markov chains, or by direct algebraic analysis of the mass function.

As with the occupancy distribution, it is also possible to derive interesting mixture results using the negative occupancy distribution, which generalise well-known mixture representations for the negative binomial distribution. We will derive the mixture characterisation directly through the mass function in this case.

**THEOREM 14 (Negative occupancy distribution is a negative binomial mixture of negative occupancy distributions):** The negative occupancy distribution satisfies the equation:

$$\text{NegOcc}(t|m,k,\gamma\theta) = \sum_{r=0}^{t} \text{NegBin}(t-r|k+r, 1-\theta) \cdot \text{NegOcc}(r|m,k,\gamma)\,.$$

In the special case where $\gamma = 1$ we obtain the useful mixture equation:

$$\text{NegOcc}(t|m,k,\theta) = \sum_{r=0}^{t} \text{NegBin}(t-r|k+r, 1-\theta) \cdot \text{NegOcc}(r|m,k)\,.$$

Theorem 14 is the negative occupancy analogue to Theorem 12 for the occupancy distribution. This theorem also has a simple intuition when we interpret it as involving two independent events resulting in "falling through" the bins. The negative binomial distribution in the mixture gives the component of the excess hitting time that is attributable to the new event, and each of these terms is multiplied by the negative occupancy distribution without the probability of that event incorporated. We have stated the theorem in a general form, but the most important case occurs when $\gamma = 1$, which allows us to write the negative occupancy distribution as a negative binomial mixture of the classical negative occupancy distribution. Again, this latter mixture is especially useful for computational purposes.

We have seen that the negative occupancy distribution generalises the negative binomial distribution, so it is useful to compare the mixture characterisation in Theorem 14 to the known characterisations of the negative binomial distribution. Taking $m \rightarrow \infty$ in Theorem 14 gives $\text{NegOcc}(t|m,k,\gamma) \rightarrow \text{NegBin}(t|k,1-\gamma)$ so that Theorem 14 reduces asymptotically down to the well-known negative binomial mixture:

$$\text{NegBin}(t|k,1-\gamma\theta) = \sum_{r=0}^{t} \text{NegBin}(t-r|k+r, 1-\theta) \cdot \text{NegBin}(r|k,1-\gamma)\,.$$

Each of the above mixtures extends known mixture characterisations for the binomial, Poisson, or negative binomial distributions. To complete our analysis of mixture characterisations, we will derive one final mixture result using the spillage distribution. This mixture does not extend any other well-known mixture results. (Taking $m \rightarrow \infty$ reduces the spillage distribution to a point-mass distribution on zero, so in this limiting case the mixture reduces to a trivial assertion that $\text{Bin}(s|n,\theta) = \text{Bin}(s|n,\theta)$.)

**THEOREM 15 (Binomial distribution is a spillage mixture of occupancy distributions):** The binomial distribution satisfies the equation:

$$\text{Bin}(s|n,\theta) = \sum_{k=0}^{s} \text{Spillage}\left(s-k \middle| n,k,m \cdot \frac{1-\theta}{\theta}\right) \cdot \text{Occ}(k|n,m,\theta) .$$

We now have a reasonably complete set of mixture results that relate the various occupancy distributions to other well-known discrete distributions. In particular, the most useful mixture results here are Theorems 12 and 14, which allow us to generate the occupancy distribution and the negative occupancy distribution as mixtures of their classical versions. (In the first case this is a binomial mixture and in the second case it is a negative binomial mixture.) These results are useful for computational purposes, since they allow us to generate the extended distributions from their classical counterparts.

## 7. Application to "coverage" analysis in bootstrapping/resampling problems

The extended occupancy problem and the three distributions discussed in this paper arise when we undertake simple-random-sampling without replacement (SRSWOR) from a finite set of objects. One statistical context in which this occurs is when we use resampling methods such as "bootstrap" estimation (see e.g., Hall 1992). In this context we may wish to examine the coverage of the data points in the original sample, which leads us to the classical occupancy problem. Consequently, various aspects of the coverage of the original sample are described by the occupancy distributions we have looked at in this paper. In particular, the marginal and conditional distributions of the coverage of the original sample are described by the extended occupancy distribution, and our other distributions describe related aspects of the problem.

The formal description of occupancy analysis in resampling problems is fairly straightforward. In order to stick with our existing notation throughout this paper, suppose we have an initial sample of data points $\boldsymbol{x} = (x_1, \dots, x_m)$ and we decide to resample $n$ data points via SRSWOR.

We note that it is usual in bootstrapping to generate resamples that are of the same size as the original sample (i.e., with $n = m$). We will proceed in greater generality because our coverage analysis applies just as well to resampling that does not impose this restriction. For simplicity, we will also assume that all the original data points are distinct, such as would occur when the underlying distribution is continuous. (Analysis can be extended to the case where there are duplicate values in the original sample, but certain aspects of the problem then go beyond the extended occupancy problem.) Formally, bootstrapping works by generating a resampled data vector $\boldsymbol{y} = (y_1, \dots, y_n)$ where the elements are:

$$y_i = x_{U_i} \qquad\qquad U_1, \dots, U_n \sim \text{IID U}\{1, \dots, m\}.$$

Let $\mathcal{J}_n \equiv \bigcup_{i=1}^{n} \{U_i\} \subseteq \{1, \dots, m\}$ be the set of data points that were resampled (as described by their indices) and let $K_n \equiv |\mathcal{J}_n|$ be the size of this set. Let $T_k \equiv \min\{t = 0,1,2, \dots | K_{k+t} = k\}$ be the number of excess resampled values required to ensure that the resampled vector includes $k$ different data points from the original sample. The set $\mathcal{J}_n$ describes the "coverage" of the original data points in the resample, and the quantity $K_n$ tells us the number of data points in the original sample that appear in the resample.

Since the indices for the resampled values are independent uniform random variables over the original indices for the data points, we can easily see that the "coverage" of the original sample is described by the classical occupancy problem. (Our notation $K_n$ and $T_k$ for the coverage of the original sample and its excess hitting time reflect the fact that these are the occupancy number and excess hitting time in the classical occupancy problem.) We can easily see that the number of resampled data points $K_n$ is an occupancy number from the classical occupancy problem, so it follows the marginal and conditional distributions:

$$\mathbb{P}(K_n = k) = \text{Occ}(k|n, m),$$

$$\mathbb{P}(K_{\acute{n}+n} = k | K_{\acute{n}} = r) = \text{Occ}(k - r | n, m - r, 1 - r/m).$$

Similarly, we can easily see that $T_k$ is an excess hitting time in the classical occupancy problem, so it follows the marginal and conditional distributions:

$$\mathbb{P}(T_k = t) = \text{NegOcc}(t|m, k),$$

$$\mathbb{P}(T_{\acute{k}+k} = t | T_{\acute{k}} = r) = \text{NegOcc}(t - r | m - r, k, 1 - r/m).$$

Since resampling follows the classical occupancy problem, in this context the "spillage" is trivial; we always have $n_{\text{eff}} = n$ so $n_{\text{eff}} - K_n = n - K_n$ with probability one.

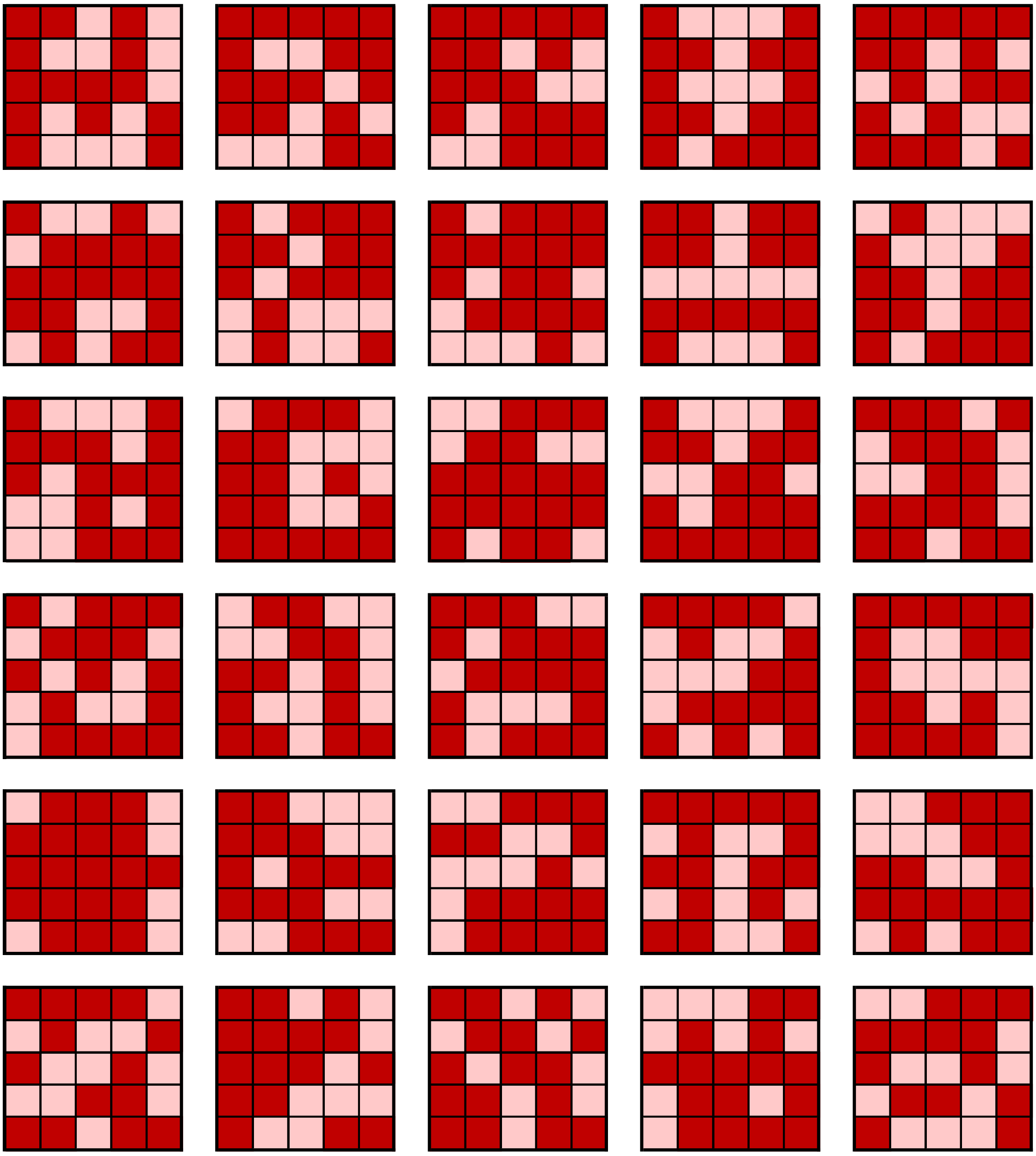

**Figure 2:** Coverage of simulated bootstrap resamples for an original sample of $n = 25$ data points using $m = 25$ resampled points. There are $S = 30$ squares showing simulated resamples. Each red square in the larger square represents a data point that is included in the resample; the occupancy number for each resample is the number of red squares.

Bootstrapping involves taking some large number of resamples from an original sample vector in order to estimate the sampling distribution of a quantity of interest. Suppose we generate bootstrap simulations $s = 1, \dots, S$ giving corresponding resample vectors $\boldsymbol{y}^{(1)}, \dots, \boldsymbol{y}^{(S)}$ using the above method. If $S$ is large then we can rely on the strong law of large numbers to assure ourselves that the empirical distributions of the coverage quantities for these resamples will converge almost surely to their true distributions:

$$\frac{1}{S}\sum_{s=1}^{S}\mathbb{I}(K_n^{(s)}=k)\overset{\text{a.s}}{\rightarrow}\text{Occ}(k|n,m),$$

$$\frac{1}{S}\sum_{s=1}^{S}\mathbb{I}(T_k^{(s)}=t)\overset{\text{a.s}}{\rightarrow}\text{NegOcc}(t|m,k).$$

In Figure 2 above we show the coverage of $S = 30$ bootstrap resamples of an original sample containing $n = 25$ data points. Each simulated resample is shown by one of the larger squares, and the occupancy number for each resample is the number of red squares in the larger square. (Since we are interested only in the coverage of the original sample here, the figure does not show how many times each value was resampled — only if it is resampled at least once or not.) As we take more and more resamples (i.e., as $S \to \infty$), the empirical distribution of the occupancy numbers will converge to the classical occupancy distribution with $n = m = 25$.

Analysis of coverage of the original sample has been used in bootstrapping analysis in order to examine issues of bias that arise from conflation of training and testing. For example, Efron and Tibshirani (1997) examine bootstrapping and cross-validation methods for estimation of error rates in binary regression modelling. In order to correct for bias arising from training and testing with the same data, they examine a special type of bootstrapping analysis formulated in Efron (1983) called the "0.632+ bootstrap method" (see also Efron 1986). Though they give their own explanation of this method, here we will examine this method in terms of our own coverage analysis. To do this, we will begin by noting that the expectation and variance of the coverage proportion are:

$$\mathbb{E}\left(\frac{K_n}{m}\right)=\left[1-\left(1-\frac{1}{n}\right)^n\right],$$

$$\mathbb{V}\left(\frac{K_n}{m}\right)=\left[(m-1)\left(1-\frac{2}{n}\right)^n+\left(1-\frac{1}{n}\right)^n-m\left(1-\frac{1}{n}\right)^{2n}\right].$$

If we take $n \to \infty$ and $m \to \infty$ subject some fixed limiting ratio $n/m \to \lambda$ then we have the asymptotic equivalence:

$$\mathbb{E}\left(\frac{K_n}{m}\right)\sim 1-e^{-\lambda} \qquad \mathbb{V}\left(\frac{K_n}{m}\right)\sim\frac{e^{-\lambda}(1-e^{-\lambda})}{m}.$$

Using the standard bootstrapping method where we use a resample that is the same size as the original sample (i.e., with $n = m$ so that $\lambda = 1$) we have convergence $K_n \to 1 - 1/e \approx 0.632$ (using the strong law of large numbers this is almost sure convergence). This is empirically evident in the resample simulations in Figure 2.

For any particular data point $i$, it can similarly be shown that $\mathbb{P}(i \in \mathcal{J}_n) \to 1 - 1/e$. Efron and Tibshirani (1997) note that error rate analysis suffers from bias when a data point used for prediction purposes (the "test point") is also included in the training sample, and this is usually dealt with in cross-validation by using the "leave one out" method (i.e., the test point is left out of the training set used for constructing its prediction). However, in bootstrapping analysis the resample will include the test point with approximate probability $0.632$. The "0.632+ bootstrap method" formulates an estimator that takes a weighted average of two estimators with known upward and downward biases, where the probabilities of inclusion/exclusion of the test point are used as weightings in the method. Our purpose here is not to recommend this method. (Indeed, the present author has a great deal of scepticism towards bootstrapping methods, but that is beyond the scope of this paper.) It is simply to note that analysis of coverage of the original data points is an important issue in the analysis of bootstrapping and resampling, and it leads to methods that take account of coverage probabilities pertaining to the original sample.

Suppose now that we look more broadly than the bootstrap, at resampling methods that may use a different number of resample values than were in the original sample (i.e., allowing for the case where $n \neq m$). One natural question in this context is how many points one should resample in order to get some stipulated minimum probability of a particular level of coverage of the original sample (e.g., including at least $k$ distinct points from the original sample). The probability that a resample of size $n$ covers at least $k$ data points in the original sample is:

$$\mathbb{P}(K_n \geq k) = \mathbb{P}(T_k \leq n - k) = \sum_{r=k}^{m} \text{Occ}(r|n, m) = \sum_{s=0}^{n-k} \text{NegOcc}(s|m, k) .$$

Given some stipulated minimum probability $0 < \phi < 1$ we can use the occupancy distributions to find the required resample size for this problem:

$$\begin{aligned} \hat{n}_k(\phi) &\equiv \min\{n \in \mathbb{N} \mid \textstyle\sum_{r=k}^{m} \text{Occ}(r|n, m) \geq \phi\} \\ &= \min\{n \in \mathbb{N} \mid \textstyle\sum_{r=k}^{m} (m)_r \cdot S(n, r) \geq \phi m^n\} . \end{aligned}$$

The value $\hat{n}_k(\phi)$ is the smallest number of resampled values required to give a probability of at least $\phi$ of an occupancy number at least $k$. Computation of this quantity allows an analyst to pre-determine the required resample size for a coverage requirement on the original sample. The special case where we seek full coverage of the original sample (i.e., $K_n = m$) is a variation of the coupon-collector problem. Investigations of this kind can be of use if an analyst wishes to undertake resampling in a manner that is likely to give some specified level of coverage of the original sample.

## 8. Summary and concluding remarks

Our goal in this paper has been to derive and discuss three interesting distributions arising from the extended occupancy problem. This problem can be framed as a pure-birth Markov chain describing the evolution of the occupancy number as more and more balls are added to a fixed number of bins, with some fixed probability of occupancy for each ball. The three distributions we have examined arise to describe the behaviour of various aspects of this Markov chain — the occupancy number, the excess hitting time for the occupancy number, and the "spillage" describing the difference between the effective number of balls and the occupancy number. It is interesting that all three distributions involve the noncentral Stirling numbers of the second kind, and the first two distributions generalise other well-known discrete distributions.

Setting aside their statistical derivation, the mathematical form of these distributions is also interesting, insofar as each distribution can be framed as the distributional analogue to a well-known equation for the noncentral Stirling numbers of the second kind (i.e., each corresponds to a normed version of a summation result involving the noncentral Stirling numbers of the second kind). The occupancy distribution arises as the distributional analogue to the equation expressing a power-sum as a sum of falling factorials of one of the values, with the noncentral Stirling numbers of the second kind arising as the coefficients. The negative occupancy distribution arises as the distributional analogue to the ordinary generating function of the noncentral Stirling numbers of the second kind. Finally, the spillage distribution arises as the distributional analogue to the equation expressing the noncentral Stirling numbers of the second kind as a weighted sum of the (central) Stirling numbers of the second kind.

The occupancy distributions in this paper arise in contexts where we undertake simple-random-sampling without replacement from a finite set of items, and we then examine the "occupancy number" and related quantities pertaining to the sample. This also arises in bootstrapping and other resampling techniques, where we can use the occupancy distributions to describe the stochastic behaviour of various quantities looking at the coverage of the original sample.

We hope that this inquiry has given the reader an appreciation for the various ways that the noncentral Stirling numbers of the second kind arise in the extended occupancy problem, and has likewise given an appreciation for the fact that this simple problem generates distributions that provide analogues to a wide range of equations involving the noncentral Stirling numbers of the second kind. Both the statistical and mathematical aspects of these three distributions are interesting, and they provide a broad class of discrete distributional forms created from the noncentral Stirling numbers of the second kind. It is particularly interesting that our first two distributions provide generalisations of the binomial and negative binomial distributions, with an additional parameter $0 \leq m \leq \infty$ that has the effect of "squashing" the occupancy number when we impose a finite value. Imposing a finite number of bins on the extended occupancy problem will tend to give a lower value of the occupancy number, and thus a higher value of the excess hitting time, than would be the case if balls were allocated among an infinite number of bins. This also leads to a non-trivial distribution for the "spillage", measuring the difference between the effective number of balls and the occupancy number. As a reference, we give some tables below that summarise our three distributions, and summarise the mixture results involving these distributions.

| Occupancy Distribution | |
|---|---|
| Mass function | $$\text{Occ}(k\|n,m,\theta) = \frac{\theta^n}{m^n} \cdot (m)_k \cdot S\left(n,k,m \cdot \frac{1-\theta}{\theta}\right)$$ |
| Comparison to binomial | $$\text{Occ}(k\|n,m,\theta) = \frac{(m)_k}{m^k} \cdot \Pi\left(n,k,m \cdot \frac{1-\theta}{\theta}\right) \times \text{Bin}(k\|n,\theta)$$ |
| Distributional analogue of… | $$(m+\phi)^n = \sum_{k=0}^{n} (m)_k \cdot S(n,k,\phi)$$ |
| Stochastic dominance | $n \leq n' \implies F(k\|n,m,\theta) \geq F(k\|n',m,\theta)$ (strict if $n < n'$ and $m > 1$)<br>$m \leq m' \implies F(k\|n,m,\theta) \geq F(k\|n,m',\theta)$ (strict if $m < m'$ and $n > 1$)<br>$\theta \leq \theta' \implies F(k\|n,m,\theta) \geq F(k\|n,m,\theta')$ (strict if $\theta < \theta'$) |
| Central moments | Central moment formulae use $E_r \equiv (1-\theta r/m)^n$. They are:<br>$$\mu_{n,m,\theta} = m(1-E_1),$$<br>$$\sigma^2_{n,m,\theta} = m[(m-1)E_2 + E_1 - mE_1^2],$$<br>$$\gamma_{n,m,\theta} = -\frac{E_1 - 3E_2 + 2E_3 + m\begin{pmatrix} 3(E_2 - E_1^2) \\ +2m(E_1^3 - E_1E_2) \\ +(m-3)(E_3 - E_1E_2) \end{pmatrix}}{m^{1/2}[(m-1)E_2 + E_1 - mE_1^2]^{3/2}},$$<br>$$\kappa_{n,m,\theta} = \frac{\begin{pmatrix} E_1 - 4mE_1^2 + 6m^2E_1^3 - 3m^3E_1^4 \\ +7(m-1)E_2 + 6(m-1)(m-2)E_3 \\ +(m-1)(m-2)(m-3)E_4 \\ -12m(m-1)E_1E_2 + 6m^2(m-1)E_1^2E_2 \\ -4m(m-1)(m-2)E_1E_3 \end{pmatrix}}{m[(m-1)E_2 + E_1 - mE_1^2]^2}.$$ |

**Table 1:** Summary of the occupancy distribution

| Negative Occupancy Distribution | |
|---|---|
| Mass function | $$\text{NegOcc}(t\|m,k,\theta) = \frac{\theta^{k+t}}{m^{k+t}} \cdot (m)_k \cdot S\left(k+t-1,k-1,m \cdot \frac{1-\theta}{\theta}\right)$$ |
| Comparison to negative binomial | $$\text{NegOcc}(t\|m,k,\theta) = \frac{(m)_k}{m^k} \cdot \Pi\left(k+t-1,k-1,m \cdot \frac{1-\theta}{\theta}\right) \times \text{NegBin}(t\|k,1-\theta)$$ |
| Distributional analogue of… | $$\sum_{n=k}^{\infty} S(n,k,\phi) \cdot x^n = \frac{x^k}{\prod_{i=0}^{k}(1-(i+\phi)x)}$$ |
| Stochastic dominance | $m \leq m' \implies F(t\|m,k,\theta) \leq F(t\|m',k,\theta)$ (strict if $m < m'$ and $k > 1$ or $\theta < 1$)<br>$k \leq k' \implies F(t\|m,k,\theta) \geq F(t\|m,k',\theta)$ (strict if $k < k'$)<br>$\theta \leq \theta' \implies F(t\|m,k,\theta) \leq F(t\|m,k,\theta')$ (strict if $\theta < \theta'$) |

**Table 2:** Summary of the negative occupancy distribution

| Spillage Distribution | |
|---|---|
| Mass function | $\text{Spillage}(r\|n,k,\phi) = \binom{n}{k+r} \cdot \phi^{n-k-r} \cdot \frac{S(k+r,k)}{S(n,k,\phi)}$ |
| Distributional analogue of… | $S(n,k,\phi) = \sum_{r=k}^{n} \binom{n}{r} \cdot \phi^{n-r} \cdot S(r,k)$ |

**Table 3:** Summary of the spillage distribution

| Mixtures involving occupancy distributions | |
|---|---|
| Occupancy as a binomial mixture of occupancy distributions | $\text{Occ}(k\|n,m,\gamma\theta) = \sum_{r=0}^{n} \text{Bin}(r\|n,\theta) \cdot \text{Occ}(k\|r,m,\gamma)$ |
| Binomial as a Poisson mixture of occupancy distributions | $\text{Bin}(k\|m,\gamma) = \sum_{r=0}^{\infty} \text{Pois}\left(r \middle\| m \cdot \frac{\|\ln(1-\gamma)\|}{\theta}\right) \cdot \text{Occ}(k\|r,m,\theta)$ |
| Binomial as a spillage mixture of occupancy distributions | $\text{Bin}(s\|n,\theta) = \sum_{k=0}^{s} \text{Spillage}\left(s-k \middle\| n,k,m \cdot \frac{1-\theta}{\theta}\right) \cdot \text{Occ}(k\|n,m,\theta)$ |
| Negative occupancy as a negative binomial mixture of negative occupancy distributions | $\text{NegOcc}(t\|m,k,\gamma\theta) = \sum_{r=0}^{t} \text{NegBin}(t-r\|k+r,1-\theta) \cdot \text{NegOcc}(r\|m,k,\gamma)$ |

**Table 4:** Summary of mixture results involving the occupancy distributions

**Declarations**

| | |
|---|---|
| **Availability of data and materials** | Not applicable |
| **Competing interests** | The author declares that he has no competing interests |
| **Funding** | There was no specific funding for this paper or project |
| **Author's contributions** | The sole author undertook all work on the present paper (i.e., conception, analysis and writing for the paper). |
| **Acknowledgements** | The author would like to thank two anonymous referees at the *Journal for Statistical Distributions and Applications* who provided useful suggestions to improve a previous draft of this paper. |

## Appendix I: Noncentral Stirling numbers of the second kind

In this appendix we set out some basic material on noncentral Stirling numbers of the second kind, which run through much of our analysis. At the outset, we note a trap for the unwary, which is that there are two alternative definitions of these numbers in the literature. Our definition of these numbers is equivalent to the one presented in Charalambides (2005, p. 75, Eqn. 2.7), using the values $S(n,k,\phi)$ that satisfy the expansion:

$$(t+\phi)^n = \sum_{k=0}^{n} S(n,k,\phi)\,(t)_k,$$

where the values $(t)_k = \prod_{i=0}^{k-1}(t-i)$ are falling factorials. We can see from this equation that the noncentral Stirling numbers of the second kind are the coefficients of the expansion that converts a power of a sum into a sum of falling factorials. The values $n$ and $k$ are non-negative integers and the value $\phi$ is a real **noncentrality parameter**. In the special case where $\phi = 0$ we obtain the standard (central) Stirling numbers of the second kind, which satisfy the simpler equation $t^n = \sum_{k=0}^{n} S(n,k)\,(t)_k$. The Stirling numbers of the second kind arise in discrete analysis, mostly through these conversion equations. Using the forward difference operator $\Delta$ these numbers can be written as $S(n,k,\phi) = [\Delta^k t^n/k!]_{t=\phi}$ (p. 76, Eqn. 2.21). This method of framing the operator is also used by Uppuluri and Carpenter (1971) (pp. 316).

**NOTE:** As stated, there are two alternative definitions of these numbers in the literature. The two alternative definitions of the noncentral Stirling numbers of the second kind correspond to two directions in which the noncentrality parameter can be expressed. The definition used here differs from the definition used in Koutras (1982, pp. 81-84), where the author measures non-centrality in the other direction, using the noncentrality parameter $a = -\phi$. When using this latter definition, all the resulting equations involving these numbers include negative signs attached to the noncentrality parameter (and powers of negative signs in the various expansions). Either definition is perfectly serviceable; they merely express noncentrality in different directions. However, when dealing with the extended occupancy problem, it is more natural to use the definition that avoids adding unnecessary negative signs into our various equations. The reader should of course be careful to ensure that they bear our definition in mind when applying any of the equations in this paper; use of our equations with the alternative definition of these numbers leads to errors. □

Our goal in this paper concerns analysis of the extended occupancy problem, so we will not derive the properties of the Stirling numbers. Charalambides (2005, pp. 73-96) provides a comprehensive introduction to these numbers, with relevant derivations of properties, and this is a useful reference for the interested reader. We give a selection of useful results here. There are several explicit forms for the noncentral Stirling numbers of the second kind, which are each useful in different contexts. One explicit form (p. 85, Eqn. 2.26) is:

$$S(n, k, \phi) = \frac{1}{k!}\sum_{i=0}^{k} \binom{k}{i} (-1)^{k-i}(i + \phi)^n .$$

The noncentral Stirling numbers of the second kind can also be written explicitly as a sum of powers of values that are increments of the noncentrality parameter (pp. 93-94), as:

$$S(n, k, \phi) = \sum_{\boldsymbol{j}} \phi^{n-k-j_1-\cdots-j_k}(\phi + 1)^{j_1} \ldots (\phi + k)^{j_k} ,$$

where the sum is taken over the set of all vectors $\boldsymbol{j} = (j_1, \ldots, j_k)$ of non-negative integer values with sum $\sum_i j_i \leq n - k$. In practice, the noncentral Stirling numbers of the second kind are usually computed via the following triangular recursive equation (p. 88, Eqn 2.34):

$$S(n + 1, k, \phi) = (k + \phi) \cdot S(n, k, \phi) + S(n, k - 1, \phi),$$

combined with the initial conditions:[13]

$$\begin{aligned} S(n, 0, \phi) &= \phi^n && \text{for all } n \geq 0, \\ S(n, n, \phi) &= 1 && \text{for all } n \geq 0, \\ S(n, k, \phi) &= 0 && \text{for all } k > n. \end{aligned}$$

Many properties of the Stirling numbers are derived using generating functions. The ordinary generating function for these numbers (p. 93, Eqn. 2.40) is:

$$\sum_{n=k}^{\infty} S(n, k, \phi) \cdot x^n = \frac{x^k}{\prod_{i=0}^{k}(1 - (i + \phi)x)} \qquad |x| < \frac{1}{\phi + k}.$$

The exponential generating function for these numbers (p. 80, Eqn. 2.24) is:

$$\sum_{n=k}^{\infty} S(n, k, \phi) \cdot \frac{x^n}{n!} = e^{\phi x} \cdot \frac{(e^x - 1)^k}{k!} \qquad x \in \mathbb{R}.$$

A related generating function (p. 80, Eqn. 2.22) is:

$$\sum_{n=0}^{\infty}\sum_{k=0}^{n} S(n, k, \phi)(t)_k \cdot \frac{x^n}{n!} = e^{(t+\phi)x}.$$

[13] The second of these conditions is a restatement of an algebraic identity given in Ruiz (1996).

Our analysis looks at recursive properties of distributions involving the noncentral Stirling numbers of the second kind. We first note that these numbers can be written in terms of the (central) Stirling numbers of the second kind via the equation (p. 77, Eqn. 2.16):

$$S(n,k,\phi) = \sum_{i=k}^{n} \binom{n}{i} \phi^{n-i} \cdot S(i,k) .$$

**Lemma 1 (Recursive and differential equations):** We have:

$$S(n+1,k,\phi) = (k+\phi) \cdot S(n,k,\phi) + S(n,k-1,\phi),$$

$$S(n+1,k,\phi) = \phi \cdot S(n,k,\phi) + S(n,k-1,\phi+1),$$

$$S(n,k,\phi+1) = (k+1) \cdot S(n,k+1,\phi) + S(n,k,\phi),$$

$$\frac{\partial}{\partial\phi} S(n,k,\phi) = n \cdot S(n-1,k,\phi).$$

**Lemma 2 (Telescoping equation):** We have:

$$S(n+1,k,\phi) = \sum_{r=0}^{n-k+1} (k+\phi)^r \cdot S(n-r,k-1,\phi) .$$

**Lemma 3 (Moving the noncentrality parameter):** For all $\phi > 0$ and $\phi' > 0$ we have:

$$S(n,k,\phi') = \sum_{r=0}^{n-k} \binom{n}{r} (\phi' - \phi)^r \cdot S(n-r,k,\phi) .$$

To facilitate our analysis we will formulate a function giving a scaled version of the noncentral Stirling numbers of the second kind in the case where $\phi > 0$. For fixed values of $n$ and $k$ the second of these equations is an explicit polynomial in $\phi$. It can be written in long-form as:

$$S(n,k,\phi) = \binom{n}{k} \cdot \phi^{n-k} + \cdots + \binom{n}{n-1} \cdot S(n-1,k) \cdot \phi + S(n,k).$$

Dividing by the leading term of this polynomial gives the **scaled Stirling function**:[14]

$$\Pi(n,k,\phi) \equiv \frac{S(n,k,\phi)}{\phi^{n-k}} \Big/ \binom{n}{k} \qquad \text{for } 0 \le n,k \le \infty \text{ and } 0 < \phi \le \infty.$$

[14] In the cases where $k > n$ we define $\Pi(n,k,\phi) \equiv 0$ by convention. We allow infinite inputs for the argument variables, and in these cases we define the output of the function by the corresponding limits.

The scaled Stirling function can be written in various useful forms, among which is:

$$\Pi(n,k,\phi) \equiv \frac{S(n,k,\phi)}{\phi^{n-k}} \bigg/ \binom{n}{k} = \sum_{i=k}^{n} \frac{\binom{n}{i}}{\binom{n}{k}} \cdot S(i,k) \cdot \phi^{k-i}$$

$$= \sum_{i=k}^{n} \frac{(n-k)!\,k!}{(n-i)!\,i!} \cdot S(i,k) \cdot \phi^{k-i}$$

$$= \sum_{i=0}^{n-k} \frac{(n-k)!\,k!}{(n-k-i)!\,(k+i)!} \cdot S(k+i,k) \cdot \phi^{-i}$$

$$= \sum_{i=0}^{n-k} \frac{(n-k)_i}{(k+i)_i} \cdot S(k+i,k) \cdot (1/\phi)^i \,.$$

(Here we use the "falling factorials" $(m)_t = \prod_{i=0}^{t-1}(m-i)$ to expand binomial coefficients.) This is a rational function of $\phi$ that reduces to a polynomial in $1/\phi$. The polynomial is scaled so that its constant term is unity. The scaled Stirling function is shown below in Figure 3 for fixed values of $n$ and $k$, plotted as a function of $\phi$ on a dual-logarithmic scale. As $\phi \to 0$ we have asymptotic equivalence $\log \Pi(n,k,\phi) \sim \log S(n,k) - \log\binom{n}{k} - (n-k)\log\phi$, which means that on the dual-logarithmic scale the curves left-converge to lines with slope $n-k$.

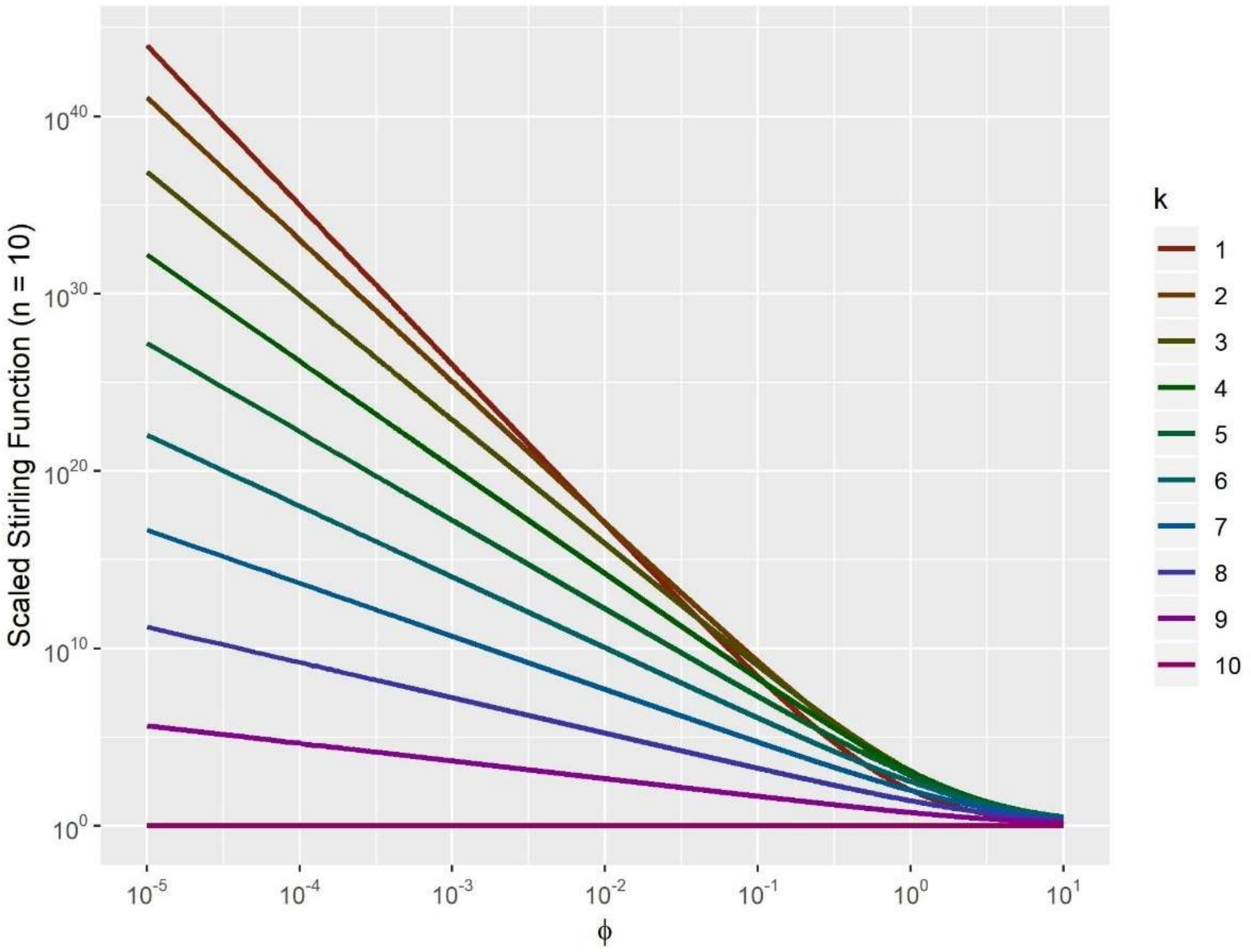

**Figure 3:** Plot of the scaled Stirling function for $n = 10$ and $k = 1, \dots, 10$.

The support of the scaled Stirling function occurs over the subdomain $0 \leq k \leq n \leq \infty$. Within this support, the function can be written in its "canonical form" as follows. Define the set $\mathcal{J}_{n,k}$ containing all vectors $\boldsymbol{j} = (j_1, \dots, j_k)$ of non-negative integer values that sum to $\sum_i j_i \leq n - k$. The number of elements in this set is equal to the number of ways to partition $n - k$ objects into $k + 1$ subsets. Using the classic "stars and bars" combinatorial method of Feller (1950) we can represent this as $k$ "bars" in the gaps between $n - k$ "stars", which yields $|\mathcal{J}_{n,k}| = \binom{n}{k}$. (Note that this set is empty if we are outside the support of $\Pi$.) Since $\phi > 0$ over the domain of our function, we can take the corresponding equation for the noncentral Stirling numbers of the second kind and divide it through by $\phi^{n-k}$ to get the canonical form:

$$\Pi(n, k, \phi) = \frac{1}{|\mathcal{J}_{n,k}|} \sum_{\boldsymbol{j} \in \mathcal{J}_{n,k}} \left(1 + \frac{1}{\phi}\right)^{j_1} \dots \left(1 + \frac{k}{\phi}\right)^{j_k}.$$

This equation can be given a probabilistic interpretation, by observing that averaging over a set of vector inputs is equivalent to finding the expected value of that function for a random vector that is uniformly distributed over the aforementioned set. To formalise this idea, we define the function $H(\mathbf{J}_k, \phi) \equiv (1 + 1/\phi)^{J_1} \dots (1 + k/\phi)^{J_k}$. Taking $\mathbf{J}_k \sim \mathrm{U}(\mathcal{J}_{n,k})$ as a random vector that is uniformly distributed over the set $\mathcal{J}_{n,k}$, we then have:

$$\Pi(n, k, \phi) = \mathbb{E}(H(\mathbf{J}_k, \phi)).$$

In order to derive certain properties of the scaled Stirling function, it is useful to decompose this function using the law of total probability, by conditioning on $\dot{J} \equiv \sum_i J_i$. To formalise this, for each $\ell = k, \dots, n$ we let $\mathcal{U}_{\ell,k}$ denote the set of all vectors $\boldsymbol{j} = (j_1, \dots, j_k)$ containing non-negative integer values with sum $\sum_i j_i = \ell - k$. Using the same "stars and bars" argument used above, we can establish that $|\mathcal{U}_{\ell,k}| = \binom{\ell-1}{k-1}$, and so we have:

$$\mathbb{P}(\dot{J} = \ell - k) = \mathbb{P}(\mathbf{J}_k \in \mathcal{U}_{\ell,k}) = \frac{|\mathcal{U}_{\ell,k}|}{|\mathcal{J}_{n,k}|} = \frac{\binom{\ell-1}{k-1}}{\binom{n}{k}}.$$

Defining the conditional expectation function $\Lambda(\ell, k, \phi) \equiv \mathbb{E}(H(\mathbf{J}_k, \phi) | \dot{J} = \ell - k)$, and using the law of total probability, we can decompose the scaled Stirling function as:

$$\begin{aligned} \Pi(n, k, \phi) = \mathbb{E}(H(\mathbf{J}_k, \phi)) &= \sum_{\ell=k}^{n} \mathbb{P}(\dot{J} = \ell - k) \cdot \mathbb{E}(H(\mathbf{J}_k, \phi) | \dot{J} = \ell - k) \\ &= \sum_{\ell=k}^{n} \frac{\binom{\ell-1}{k-1}}{\binom{n}{k}} \cdot \Lambda(\ell, k, \phi). \end{aligned}$$

We can now establish some useful asymptotic properties for the scaled Stirling function, which we state in conjunction with corresponding limiting results.

**Lemma 4 (Monotonicity, concavity, and asymptotics):** The scaled Stirling function has the following monotonicity, concavity, and asymptotic properties over its support.

(a) If $k < n$ the function is monotonically decreasing in $\phi$ with limits:

$$\lim_{\phi \downarrow 0} \Pi(n, k, \phi) = \infty \qquad \lim_{\phi \uparrow \infty} \Pi(n, k, \phi) = 1.$$

(b) The function is monotonically increasing in $n$ with limits:

$$\Pi(0, k, \phi) = 1 \qquad \lim_{n \uparrow \infty} \Pi(n, k, \phi) = \infty.$$

(c) The function has limits in $k$ given by:

$$\Pi(n, 0, \phi) = 1 \qquad \Pi(n, n, \phi) = 1.$$

## Appendix II: Proof of Theorems

**PROOF OF LEMMA 1:** To establish the first and second equations we use the expansion:

$$S(n+1,k,\phi) = \frac{1}{k!}\sum_{i=0}^{k}\binom{k}{i}(-1)^{k-i}\cdot(i+\phi)^{n+1}.$$

For the first equation in the lemma we extract the term $(i+\phi) = (k+\phi) - (k-i)$ to obtain the two parts of required sum, and for the second we extract the term $(i+\phi) = (\phi) + (i)$ to obtain the two parts of the required sum. Equating the right-hand-sides of the first and second equations and simplifying gives $S(n,k-1,\phi+1) = k\cdot S(n,k,\phi) + S(n,k-1,\phi)$, which gives the third equation. Finally, differentiating with respect to $\phi$ gives:

$$\begin{aligned}\frac{\partial}{\partial\phi}S(n,k,\phi) &= \frac{1}{k!}\sum_{i=0}^{k}\binom{k}{i}(-1)^{k-i}\cdot\frac{\partial}{\partial\phi}(i+\phi)^{n} \\ &= n\cdot\frac{1}{k!}\sum_{i=0}^{k}\binom{k}{i}(-1)^{k-i}\cdot(i+\phi)^{n-1} \\ &= n\cdot S(n-1,k,\phi),\end{aligned}$$

which establishes the last equation. ■

**PROOF OF LEMMA 2:** We will prove the theorem by induction. For $n = k-1$ the equation reduces to $S(k,k,\phi) = S(k,k,\phi)$ and for $n < k-1$ the equation reduces to $0 = 0$. Assuming the equation holds for some $n$ we can apply the first equation of Lemma 1 to obtain:

$$\begin{aligned}S(n+2,k,\phi) &= (k+\phi)\cdot S(n+1,k,\phi) + S(n+1,k-1,\phi) \\ &= \sum_{r=0}^{n-k+1}(k+\phi)^{r+1}\cdot S(n-r,k-1,\phi) + S(n+1,k-1,\phi) \\ &= \sum_{r=1}^{n-k+2}(k+\phi)^{r}\cdot S(n+1-r,k-1,\phi) + S(n+1,k-1,\phi) \\ &= \sum_{r=0}^{n-k+2}(k+\phi)^{r}\cdot S(n+1-r,k-1,\phi).\end{aligned}$$

This establishes the inductive step, which completes the proof. ■

**PROOF OF LEMMA 3:** Using the expansion of the noncentral Stirling numbers of the second kind in terms of the (central) Stirling numbers of the second kind, we have:

$$S(n,k,\phi') = \sum_{i=k}^{n} \binom{n}{i} \phi'^{n-i} \cdot S(i,k).$$

Taking $\phi' = (\phi' - \phi) + \phi$ and expanding via the binomial theorem then gives:

$$\begin{aligned}
S(n,k,\phi') &= \sum_{i=k}^{n} \binom{n}{i} (\phi' - \phi + \phi)^{n-i} \cdot S(i,k) \\
&= \sum_{i=k}^{n} \binom{n}{i} \left( \sum_{r=0}^{n-i} \binom{n-i}{r} (\phi' - \phi)^r \phi^{n-i-r} \right) \cdot S(i,k) \\
&= \sum_{i=k}^{n} \sum_{r=0}^{n-i} \binom{n-i}{r} \binom{n}{i} (\phi' - \phi)^r \phi^{n-i-r} \cdot S(i,k) \\
&= \sum_{r=0}^{n-k} \sum_{i=k}^{n-r} \binom{n-r}{i} \binom{n}{r} (\phi' - \phi)^r \phi^{n-i-r} \cdot S(i,k) \\
&= \sum_{r=0}^{n-k} \binom{n}{r} (\phi' - \phi)^r \sum_{i=k}^{n-r} \binom{n-r}{i} \phi^{n-r-i} \cdot S(i,k) \\
&= \sum_{r=0}^{n-k} \binom{n}{r} (\phi' - \phi)^r \cdot S(n-r,k,\phi),
\end{aligned}$$

which was to be shown. ■

**PROOF OF LEMMA 4:** We establish the monotonicity, concavity, and asymptotic properties of the function in the specified order, using the canonical form of the function.

(a) Each of the terms $(1 + i/\phi)$ for $i = 1, \ldots, k$ are clearly monotonically decreasing in $\phi$. Since each of the indices $j_1, \ldots, j_k$ are non-negative (and at least some of these elements will be strictly positive over at least some of the elements of the sum) the entire sum is therefore monotonically decreasing in $\phi$, which establishes the monotonicity requirement in the theorem. For the lower limit we have:

$$\begin{aligned}
\lim_{\phi \to 0} \Pi(n,k,\phi) &= \lim_{\phi \to 0} \sum_{\boldsymbol{j}} \left(1 + \frac{1}{\phi}\right)^{j_1} \cdots \left(1 + \frac{k}{\phi}\right)^{j_k} \Big/ \binom{n}{k} \\
&\geq \lim_{\phi \to 0} \sum_{\boldsymbol{j}} \left(1 + \frac{1}{\phi}\right)^{j_1 + \cdots + j_k} \Big/ \binom{n}{k}
\end{aligned}$$

$$\geq \lim_{\phi \to 0} \left(1 + \frac{1}{\phi}\right)^{n-k} \Big/ \binom{n}{k} = \infty.$$

Upper limit is:

$$\lim_{\phi \to \infty} \Pi(n, k, \phi) = \lim_{\phi \to \infty} \sum_{\boldsymbol{j}} \left(1 + \frac{1}{\phi}\right)^{j_1} \cdots \left(1 + \frac{k}{\phi}\right)^{j_k} \Big/ \binom{n}{k}$$

$$= \sum_{\boldsymbol{j}} 1^{j_1} \cdots 1^{j_k} / \binom{n}{k}$$

$$= \binom{n}{k} / \binom{n}{k} = 1.$$

(b) Consider a problem with a fixed value of $k$ but variable value of $n$. Suppose we generate a sequence of independent random vectors $\mathbf{J}_1, \mathbf{J}_2, \mathbf{J}_3, \dots$ as $\mathbf{J}_\ell \sim \text{U}(\mathcal{U}_{\ell,k})$ and another sequence of independent random variables $U_1, U_2, U_3, \dots \sim \text{IID U}\{1, \dots, k\}$ (all mutually independent). We clearly have the distributional equivalence:

$$H(\mathbf{J}_\ell, \phi) \cdot (1 + U_\ell/\phi) \overset{\text{Dist}}{\sim} H(\mathbf{J}_{\ell+1}, \phi).$$

We also have:

$$\mathbb{E}(1 + U_\ell/\phi) = \frac{1}{k} \sum_{i=1}^{k} \left(1 + \frac{i}{\phi}\right) = 1 + \frac{1}{\phi} \cdot \frac{k+1}{2}.$$

It follows that:

$$\begin{aligned}
\Lambda(\ell + 1, k, \phi) &= \mathbb{E}(H(\mathbf{J}, \phi) | \dot{J} = \ell - k + 1) \\
&= \mathbb{E}\big(H(\mathbf{J}_{\ell+1}, \phi)\big) \\
&= \mathbb{E}\big(H(\mathbf{J}_\ell, \phi) \cdot (1 + U_\ell/\phi)\big) \\
&= \mathbb{E}(1 + U_\ell/\phi) \cdot \mathbb{E}\big(H(\mathbf{J}_\ell, \phi)\big) \\
&= \mathbb{E}(1 + U_\ell/\phi) \cdot \mathbb{E}(H(\mathbf{J}, \phi) | \dot{J} = \ell - k) \\
&= \left(1 + \frac{1}{\phi} \cdot \frac{k+1}{2}\right) \cdot \Lambda(\ell, k, \phi).
\end{aligned}$$

Applying this to the decomposition of the scaled Stirling function gives:

$$\begin{aligned}
\Pi(n + 1, k, \phi) &= \sum_{\ell=k}^{n+1} \frac{\binom{\ell-1}{k-1}}{\binom{n}{k}} \cdot \Lambda(\ell, k, \phi) \\
&= \frac{1}{\binom{n}{k}} \cdot \Lambda(k, k, \phi) + \sum_{\ell=k+1}^{n+1} \frac{\binom{\ell-1}{k-1}}{\binom{n}{k}} \cdot \Lambda(\ell, k, \phi) \\
&= \frac{1}{\binom{n}{k}} \cdot \Lambda(k, k, \phi) + \sum_{\ell=k}^{n} \frac{\binom{\ell}{k-1}}{\binom{n}{k}} \cdot \Lambda(\ell + 1, k, \phi)
\end{aligned}$$

$$= \frac{1}{\binom{n}{k}} \cdot \Lambda(k,k,\phi) + \left(1 + \frac{1}{\phi} \cdot \frac{k+1}{2}\right) \sum_{\ell=k}^{n} \frac{\binom{\ell}{k-1}}{\binom{n}{k}} \cdot \Lambda(\ell,k,\phi)$$

$$> \frac{1}{\binom{n}{k}} \cdot \Lambda(k,k,\phi) + \left(1 + \frac{1}{\phi} \cdot \frac{k+1}{2}\right) \sum_{\ell=k}^{n} \frac{\binom{\ell-1}{k-1}}{\binom{n}{k}} \cdot \Lambda(\ell,k,\phi)$$

$$= \frac{1}{\binom{n}{k}} \cdot \Lambda(k,k,\phi) + \left(1 + \frac{1}{\phi} \cdot \frac{k+1}{2}\right) \Pi(n,k,\phi)$$

$$> \left(1 + \frac{1}{2\phi}\right) \cdot \Pi(n,k,\phi) > \Pi(n,k,\phi).$$

This establishes the monotonicity requirement in the theorem. For the lower limit we have $\Pi(k,k,\phi) = 1$ from the initial conditions for the noncentral Stirling numbers of the second kind. To obtain the upper limit we first note from the above inequalities that:

$$\Pi(n+1,k,\phi) = \prod_{\ell=k}^{n-1} \frac{\Pi(\ell+1,k,\phi)}{\Pi(\ell,k,\phi)} > \left(1 + \frac{1}{2\phi}\right)^{n-k}.$$

This inequality gives $\lim_{n\to\infty} \Pi(n,k,\phi) \geq \lim_{n\to\infty} (1 + 1/2\phi)^{n-k} = \infty$, which gives us the upper limit in the theorem.

(c) The limits in this part follow trivially from substitution into the function and use of the initial conditions for the noncentral Stirling numbers of the second kind.

**PROOF OF THEOREM 1:** Let $\mathcal{S}_n \equiv \mathcal{S}_n(\boldsymbol{x}_n)$ be the set of bins occupied by the first $n$ balls, and note that this is a function of $\boldsymbol{x}_n$. Since $|\mathcal{S}_n| = K_n = t$ we have $\mathbb{P}(X_{n+1} \notin \mathcal{S}_n | \boldsymbol{x}_n) = 1 - t/m$, and we also have:

$$K_{n+1} - K_n = \mathbb{I}(X_{n+1} \notin \mathcal{S}_n, X_{n+1} \neq \bullet).$$

Hence, we have:

$$\begin{aligned} \mathbb{P}(K_{n+1} - t = 1 | \boldsymbol{x}_n) &= \mathbb{P}(X_{n+1} \notin \mathcal{S}_n, X_{n+1} \neq \bullet | \boldsymbol{x}_n) \\ &= \mathbb{P}(X_{n+1} \neq \bullet) \cdot \mathbb{P}(X_{n+1} \notin \mathcal{S}_n | \boldsymbol{x}_n) \\ &= \theta \cdot (1 - t/m). \end{aligned}$$

Since $K_{n+1} - K_n$ must be either zero or one, we then have:

$$\mathbb{P}(K_{n+1} = k + t | \boldsymbol{x}_n) = \begin{cases} 1 - \theta(1 - t/m) & \text{for } k = 0, \\ \theta(1 - t/m) & \text{for } k = 1, \\ 0 & \text{otherwise.} \end{cases}$$

which was to be shown. (Note that this function depends on $\boldsymbol{x}_n$ only through $K_n = t$ so this also justifies the intermediate statement that replaces the conditioning value.) ■

**PROOF OF THEOREM 2:** Since **P** is an upper triangle matrix, its eigenvalues are its diagonal elements, which establishes the eigenvalue matrix in the theorem. Let $\mathbf{v}_k$ be the eigenvector corresponding to the eigenvalue $\lambda_k$. To confirm that **v** is an eigenvector matrix, we need to establish that for each $k = 0,1,\dots,m$ we have the characteristic equation:

$$(\mathbf{P} - \lambda_k \boldsymbol{I})\mathbf{v}_k = \mathbf{0}.$$

Since $m\lambda_k = m - (m-k)\theta$ we have:

$$\mathbf{P} - \lambda_k \boldsymbol{I} = \frac{1}{m}\begin{bmatrix} (0-k)\theta & m\theta & 0 & \cdots & 0 & 0 \\ 0 & (1-k)\theta & (m-1)\theta & \cdots & 0 & 0 \\ 0 & 0 & (2-k)\theta & \cdots & 0 & 0 \\ \vdots & \vdots & \vdots & \ddots & \vdots & \vdots \\ 0 & 0 & 0 & \cdots & (m-k-1)\theta & \theta \\ 0 & 0 & 0 & \cdots & 0 & (m-k)\theta \end{bmatrix}.$$

Hence, the characteristic equation corresponds to the scalar equations:

$$(k-i)v_{i,k} = (m-i)v_{i+1,k} \qquad i = 0,1,\dots,m.$$

It is easily shown that substitution of the values in the theorem satisfies these equations, which establishes the stated eigenvector matrix. Now, to establish the inverse eigenvector matrix we need to establish that $\mathbf{wv} = \boldsymbol{I}$. For all $k = 0,1,\dots,m$ we have:

$$\begin{aligned}(\mathbf{wv})_{i,j} = \sum_{k=0}^{m} w_{i,k} \cdot v_{k,j} &= \sum_{k=0}^{m} \binom{m-i}{k-i} \cdot (-1)^{k-j} \binom{m-k}{j-k} \\ &= \binom{m-i}{j-i} \sum_{k=0}^{m} (-1)^{k-j} \binom{j-i}{j-k} \\ &= \binom{m-i}{j-i} \sum_{k=0}^{j} (-1)^{k-j} \binom{j-i}{j-k} \\ &= \binom{m-i}{j-i} \cdot \mathbb{I}(i=j) = \mathbb{I}(i=j).\end{aligned}$$

(In the step to the second line we use the fact that $\binom{m-i}{k-i} \cdot \binom{m-k}{j-k} = \binom{m-i}{j-i} \cdot \binom{j-i}{j-k}$, which is easily established by expanding these terms out as ratios of factorials.) This shows that the inverse eigenvalue matrix is indeed the inverse of the eigenvalue matrix. The only remaining part of the theorem is the fact that the eigenvectors are linearly independent, which follows trivially from the fact that **v** has non-zero elements on its main diagonal and upper triangle, and zero elements on its lower triangle. ■

**Proof of Theorem 3:** With a bit of algebra it can be shown that:

$$(m-k)_r \binom{m}{k}\binom{k}{i} = (m)_r \binom{m-r}{i}\binom{m-r-i}{k-i}.$$

We therefore have:

$$\begin{aligned}
\mathbb{E}((m-K_n)_r) &= \sum_{k=0}^{\infty} (m-k)_r \binom{m}{k} \sum_{i=0}^{k} \binom{k}{i} (-1)^{k-i} \cdot \left(1-\theta\cdot\frac{m-i}{m}\right)^n \\
&= \sum_{k=0}^{\infty}\sum_{i=0}^{k} (m)_r \binom{m-r}{i}\binom{m-r-i}{k-i} (-1)^{k-i} \cdot \left(1-\theta\cdot\frac{m-i}{m}\right)^n \\
&= \sum_{i=0}^{\infty} (m)_r \binom{m-r}{i} \left(1-\theta\cdot\frac{m-i}{m}\right)^n \sum_{k=i}^{\infty} \binom{m-r-i}{k-i} (-1)^{k-i} \\
&= \sum_{i=0}^{\infty} (m)_r \binom{m-r}{i} \left(1-\theta\cdot\frac{m-i}{m}\right)^n \sum_{k=0}^{m-r-i} \binom{m-r-i}{k} (-1)^{k} \\
&= \sum_{i=0}^{\infty} (m)_r \binom{m-r}{i} \left(1-\theta\cdot\frac{m-i}{m}\right)^n \cdot \mathbb{I}(i=m-r) \\
&= (m)_r \left(1-\frac{r\theta}{m}\right)^n,
\end{aligned}$$

which was to be shown. ■

**Proof of Theorem 4a:** Since $\theta > 0$ we have:

$$\frac{1-\theta}{\theta} > 0.$$

Thus, the limit $m \to \infty$ implies that $\phi \equiv m \cdot (1-\theta)/\theta \to \infty$. Applying Lemma 4 then gives:

$$\lim_{m\to\infty} \frac{(m)_k}{m^k} \cdot \Pi(n,k,\phi) = \lim_{m\to\infty} \frac{(m)_k}{m^k} \times \lim_{\phi\to\infty} \Pi(n,k,\phi) = 1 \times 1 = 1.$$

The result follows using the alternative form for the occupancy mass function (written in terms of the scaled Stirling function). ■

**Proof of Theorem 4b:** Taking the limits $n \to \infty$ and $\theta \to 0$ such that $n\theta \to \lambda$ gives:

$$\left(1-\theta\cdot\frac{m-i}{m}\right)^n = \left(1-\frac{m-i}{m}\cdot\frac{n\theta}{n}\right)^n \to \exp\left(-\frac{m-i}{m}\cdot\lambda\right).$$

Applying this limit to each term in the summation in occupancy distribution then gives:

$$\lim_{\substack{n\to\infty,\theta\to 0 \\ n\theta\to\lambda}} \text{Occ}(k|n,m,\theta) = \lim_{\substack{n\to\infty,\theta\to 0 \\ n\theta\to\lambda}} \frac{\theta^n}{m^n} \cdot (m)_k \cdot S\left(n,k,m\cdot\frac{1-\theta}{\theta}\right)$$

$$= \lim_{\substack{n\to\infty,\theta\to 0\\ n\theta\to\lambda}} \binom{m}{k} \sum_{i=0}^{k} \binom{k}{i} (-1)^{k-i} \cdot \left(1 - \theta \cdot \frac{m-i}{m}\right)^{n}$$

$$= \binom{m}{k} \sum_{i=0}^{k} \binom{k}{i} (-1)^{k-i} \cdot \exp\left(-\frac{m-i}{m} \cdot \lambda\right)$$

$$= \exp(-\lambda) \binom{m}{k} \sum_{i=0}^{k} \binom{k}{i} (-1)^{k-i} \cdot \exp\left(\frac{\lambda}{m} \cdot i\right)$$

$$= \exp(-\lambda) \binom{m}{k} \sum_{i=0}^{k} \binom{k}{i} (-1)^{k-i} \cdot \exp(\lambda/m)^{i}$$

$$= (-1)^{k} \exp(-\lambda) \binom{m}{k} \sum_{i=0}^{k} \binom{k}{i} \cdot [-\exp(\lambda/m)]^{i}$$

$$= (-1)^{k} \exp(-\lambda) \binom{m}{k} (1 - \exp(\lambda/m))^{k}$$

$$= \exp(-\lambda) \binom{m}{k} (\exp(\lambda/m) - 1)^{k}$$

$$= \exp(-\lambda/m)^{m} \binom{m}{k} (\exp(\lambda/m) - 1)^{k}$$

$$= \binom{m}{k} (1 - \exp(-\lambda/m))^{k} \exp(-\lambda/m)^{m-k}$$

$$= \mathrm{Bin}\left(k \middle| m, 1 - \exp\left(-\frac{\lambda}{m}\right)\right),$$

which was to be shown. ■

**Proof of Theorem 5a:** Using the noncentrality parameter $\phi = m \cdot (1-\theta)/\theta$ we have:

$$\frac{\theta(k+\phi)}{m} = \frac{\theta}{m}\left(k + m \cdot \frac{1-\theta}{\theta}\right) = \frac{\theta k}{m} + 1 - \theta = 1 - \theta \cdot \frac{m-k}{m}.$$

We therefore have:

$$\begin{aligned}
\mathrm{Occ}(k|n+1, m, \theta) &= \frac{\theta^{n+1}}{m^{n+1}} \cdot (m)_k \cdot S(n+1, k, \phi) \\
&= \frac{\theta}{m} \cdot \frac{\theta^{n}}{m^{n}} \cdot (m)_k \cdot [(k+\phi) \cdot S(n,k,\phi) + S(n, k-1, \phi)] \\
&= \frac{\theta(k+\phi)}{m} \cdot \mathrm{Occ}(k|n,m,\theta) + \frac{\theta(m-k+1)}{m} \cdot \mathrm{Occ}(k-1|n,m,\theta) \\
&= \left(1 - \theta \cdot \frac{m-k}{m}\right) \cdot \mathrm{Occ}(k|n,m,\theta) + \theta \cdot \frac{m-k+1}{m} \cdot \mathrm{Occ}(k-1|n,m,\theta),
\end{aligned}$$

which was to be shown. ■

**LEMMA 5:** The classical occupancy distribution satisfies the following recursive equation:

$$\text{Occ}(k|n, m+1) = \left(\frac{m}{m+1}\right)^n \cdot \frac{m+1}{m-k+1} \cdot \text{Occ}(k|n, m).$$

**PROOF OF LEMMA 5:** The standard form for the classical occupancy distribution is written in terms of the falling factorials and the Stirling numbers, as:

$$\text{Occ}(k|n, m) = \frac{(m)_k \cdot S(n, k)}{(m)^n}.$$

We therefore have:

$$\begin{aligned}
\text{Occ}(k|n, m+1) &= \frac{(m+1)_k \cdot S(n, k)}{(m+1)^n} \\
&= \left(\frac{m}{m+1}\right)^n \cdot \frac{m+1}{m-k+1} \cdot \frac{(m)_k \cdot S(n, k)}{m^n} \\
&= \left(\frac{m}{m+1}\right)^n \cdot \frac{m+1}{m-k+1} \cdot \text{Occ}(k|n, m),
\end{aligned}$$

which was to be shown. ■

**PROOF OF THEOREM 5B:** We first note that with a bit of algebra it can be shown that:

$$\phi^r \cdot \text{Bin}(r|n, \theta) = (1 - \theta + \theta\phi)^n \cdot \text{Bin}\left(r \middle| n, \frac{\theta\phi}{1 - \theta + \theta\phi}\right).$$

Thus, we have:

$$\left(\frac{m}{m+1}\right)^r \cdot \text{Bin}(r|n, \theta) = \left(1 - \frac{\theta}{m+1}\right)^n \cdot \text{Bin}\left(r \middle| n, \frac{m\theta}{1 - \theta + m}\right).$$

Applying the binomial mixture in Theorem 13 and using Lemma 5 we have:

$$\begin{aligned}
\text{Occ}(k|n, m+1, \theta) &= \sum_{r=0}^{n} \text{Bin}(r|n, \theta) \cdot \text{Occ}(k|r, m+1) \\
&= \frac{m+1}{m-k+1} \cdot \sum_{r=0}^{n} \left(\frac{m}{m+1}\right)^r \cdot \text{Bin}(r|n, \theta) \cdot \text{Occ}(k|r, m) \\
&= \frac{m+1}{m-k+1} \cdot \left(1 - \frac{\theta}{m+1}\right)^n \sum_{r=0}^{n} \text{Bin}\left(r \middle| n, \frac{m\theta}{1 - \theta + m}\right) \cdot \text{Occ}(k|r, m) \\
&= \frac{m+1}{m-k+1} \cdot \left(1 - \frac{\theta}{m+1}\right)^n \cdot \text{Occ}\left(k \middle| n, m, \frac{m\theta}{1 - \theta + m}\right),
\end{aligned}$$

which was to be shown. ■

**Proof of Theorem 5c:** Writing the occupancy distribution in explicit form we have:

$$\mathrm{Occ}(k|n,m,\theta) = \binom{m}{k}\sum_{i=0}^{k}\binom{k}{i}(-1)^{i-k}\cdot\left(1-\theta\cdot\frac{m-i}{m}\right)^{n}.$$

Differentiating the occupancy distribution with respect to $\theta$ gives:

$$\begin{aligned}
\frac{\partial}{\partial\theta}\mathrm{Occ}(k|n,m,\theta) &= \binom{m}{k}\sum_{i=0}^{k}\binom{k}{i}(-1)^{i-k}\cdot\frac{\partial}{\partial\theta}\left(1-\theta\cdot\frac{m-i}{m}\right)^{n} \\
&= -\binom{m}{k}\sum_{i=0}^{k}\binom{k}{i}(-1)^{i-k}\cdot\frac{m-i}{m}\left(1-\theta\cdot\frac{m-i}{m}\right)^{n-1}
\end{aligned}$$

$$= -\mathrm{Occ}(k|n-1,m,\theta) + \frac{k}{m}\binom{m}{k}\sum_{i=1}^{k}\binom{k-1}{i-1}(-1)^{i-k}\left(1-\theta\cdot\frac{m-i}{m}\right)^{n-1}$$

$$= -\mathrm{Occ}(k|n-1,m,\theta) + \frac{k}{m}\binom{m}{k}\sum_{i=0}^{k}\left[\binom{k}{i}-\binom{k-1}{i}\right](-1)^{i-k}\left(1-\theta\cdot\frac{m-i}{m}\right)^{n-1}$$

$$= -\frac{m-k}{m}\cdot\mathrm{Occ}(k|n-1,m,\theta) - \frac{k}{m}\binom{m}{k}\sum_{i=0}^{k-1}\binom{k-1}{i}(-1)^{i-k}\left(1-\theta\cdot\frac{m-i}{m}\right)^{n-1}$$

$$= -\frac{m-k}{m}\cdot\mathrm{Occ}(k|n-1,m,\theta)$$

$$+\frac{m-k+1}{m}\cdot\binom{m}{k-1}\sum_{i=0}^{k-1}\binom{k-1}{i}(-1)^{i-k+1}\left(1-\theta\cdot\frac{m-i}{m}\right)^{n-1}$$

$$= -\frac{m-k}{m}\cdot\mathrm{Occ}(k|n-1,m,\theta) + \frac{m-k+1}{m}\cdot\mathrm{Occ}(k-1|n-1,m,\theta),$$

which was to be shown. ■

**Proof of Theorem 6a:** Applying the first recursive equation in Theorem 5 gives:

$$\begin{aligned}
F(k|n+1,m,\theta) &= \sum_{r=0}^{k}\mathrm{Occ}(r|n+1,m,\theta) \\
&= \sum_{r=1}^{k}\theta\cdot\frac{m-r+1}{m}\cdot\mathrm{Occ}(r-1|n,m,\theta) + \sum_{r=0}^{k}\left(1-\theta\cdot\frac{m-r}{m}\right)\cdot\mathrm{Occ}(r|n,m,\theta) \\
&= \sum_{r=0}^{k-1}\theta\cdot\frac{m-r}{m}\cdot\mathrm{Occ}(r|n,m,\theta) + \sum_{r=0}^{k}\left(1-\theta\cdot\frac{m-r}{m}\right)\cdot\mathrm{Occ}(r|n,m,\theta) \\
&= \sum_{r=0}^{k}\mathrm{Occ}(r|n,m,\theta) - \theta\cdot\frac{m-k}{m}\cdot\mathrm{Occ}(k|n,m,\theta)
\end{aligned}$$

$$= F(k|n,m,\theta) - \theta \cdot \frac{m-k}{m} \cdot \text{Occ}(k|n,m,\theta).$$

Thus, for all $n \leq n'$ we have:

$$F(k|n',m,\theta) = F(k|n,m,\theta) + \sum_{i=0}^{n'-n-1} \big(F(k|n+i+1,m,\theta) - F(k|n+i,m,\theta)\big)$$

$$= F(k|n,m,\theta) - \theta \cdot \frac{m-k}{m} \cdot \sum_{i=0}^{n'-n-1} \text{Occ}(k|n+i,m,\theta) \leq F(k|n,m,\theta),$$

where the inequality is strict if $n < n'$ and $m > 1$ (and noting that $\theta > 0$). ■

**PROOF OF THEOREM 6B:** As a preliminary step, we note that the second recursive equation in Theorem 5 can be applied for the case $\theta = 1$ to give the recursive formula:

$$\text{Occ}(k|n,m+1) = \frac{m+1}{m-k+1} \cdot \left(\frac{m}{m+1}\right)^n \cdot \text{Occ}(k|n,m).$$

To facilitate our analysis, we define the function $\Psi$ by:

$$\Psi(m,r) \equiv \frac{m+1}{m-r+1} \cdot \left(\frac{m}{m+1}\right)^r \qquad \text{for all } r = 0,1,2,\dots,m.$$

It is simple to establish that $\Psi(m,0) = \Psi(m,1) = 1$ and the function is strictly decreasing after this. For all $r = 0,1,2,\dots,m-1$ we have:

$$\frac{\Psi(m,r+1)}{\Psi(m,r)} = \frac{m-r+1}{m-r} \cdot \frac{m}{m+1}.$$

Now, using the binomial mixture for the extended occupancy distribution in Theorem 12, and then applying this recursive formula, gives:

$$F(k|n,m+1,\theta) = \sum_{r=0}^{k} \text{Occ}(r|n,m+1,\theta)$$

$$= \sum_{r=0}^{k} \sum_{s=r}^{n} \text{Occ}(r|s,m+1) \cdot \text{Bin}(s|n,\theta)$$

$$= \sum_{r=0}^{k} \sum_{s=r}^{n} \text{Occ}(r|s,m+1) \cdot \text{Bin}(s|n,\theta)$$

$$= \sum_{r=0}^{k} \sum_{s=r}^{n} \frac{m+1}{m-r+1} \cdot \left(\frac{m}{m+1}\right)^s \cdot \text{Occ}(r|s,m) \cdot \text{Bin}(s|n,\theta)$$

$$\leq \sum_{r=0}^{k} \Psi(m,r) \sum_{s=r}^{n} \text{Occ}(r|s,m) \cdot \text{Bin}(s|n,\theta)$$

$$= \sum_{r=0}^{k} \Psi(m,r) \cdot \mathrm{Occ}(r|n,m,\theta) \leq \sum_{r=0}^{k} \mathrm{Occ}(r|n,m,\theta) = F(k|n,m,\theta).$$

Thus, for all $m \leq m'$ we have:

$$F(k|n,m',\theta) = F(k|n,m,\theta) + \sum_{i=0}^{m'-m-1} \big(F(k|n,m+i+1,\theta) - F(k|n,m+i,\theta)\big)$$
$$\leq F(k|n,m,\theta),$$

where the inequality is strict if $m < m'$ and $n > 1$ (and noting that $\theta > 0$). ■

**PROOF OF THEOREM 6C:** Applying the differential equation in Theorem 5 gives:

$$\frac{\partial F}{\partial \theta}(k|n,m,\theta') = \sum_{r=0}^{k} \frac{\partial}{\partial \theta} \mathrm{Occ}(r|n,m,\theta)$$
$$= -\sum_{r=0}^{k} \frac{m-r}{m} \cdot \mathrm{Occ}(r|n-1,m,\theta) + \sum_{r=1}^{k} \frac{m-r+1}{m} \cdot \mathrm{Occ}(r-1|n-1,m,\theta)$$
$$= -\sum_{r=0}^{k} \frac{m-r}{m} \cdot \mathrm{Occ}(r|n-1,m,\theta) + \sum_{r=0}^{k-1} \frac{m-r}{m} \cdot \mathrm{Occ}(r|n-1,m,\theta)$$
$$= -\frac{m-k}{m} \cdot \mathrm{Occ}(k|n-1,m,\theta).$$

Thus, for all $\theta \leq \theta'$ we have:

$$F(k|n,m,\theta') = F(k|n,m,\theta) + \int_{\theta}^{\theta'} \frac{\partial F}{\partial \theta}(k,n,m,t)dt$$
$$= F(k|n,m,\theta) - \frac{m-k}{m} \int_{\theta}^{\theta'} \mathrm{Occ}(k|n-1,m,t)dt < F(k|n,m,\theta),$$

where the inequality is strict if $\theta < \theta'$.

**PROOF OF THEOREM 7:** Proof is analogous to Theorem 4(a), using the alternative form for the negative occupancy mass function (written in terms of the scaled Stirling function). ■

**PROOF OF THEOREM 8A:** Using the second recursive equation in Theorem 5 we obtain:

$$\mathrm{NegOcc}(t|m+1,k,\theta) = \theta \cdot \frac{m-k+2}{m+1} \cdot \mathrm{Occ}(k-1|k+t-1,m+1,\theta)$$
$$= \theta \cdot \left(1 - \frac{\theta}{m+1}\right)^{k+t-1} \cdot \mathrm{Occ}\left(k-1\middle|k+t-1,m,\frac{m\theta}{1-\theta+m}\right)$$

$$= \frac{1-\theta+m}{m-k+1} \cdot \left(1 - \frac{\theta}{m+1}\right)^{k+t-1}$$

$$\times \theta \cdot \frac{m-k+1}{1-\theta+m} \cdot \text{Occ}\left(k-1 \middle| k+t-1, m, \frac{m\theta}{1-\theta+m}\right)$$

$$= \frac{1-\theta+m}{m-k+1} \cdot \left(1 - \frac{\theta}{m+1}\right)^{k+t-1} \cdot \text{NegOcc}\left(t \middle| m, k, \frac{m\theta}{1-\theta+m}\right)$$

$$= \frac{1-\theta+m}{m+1} \cdot \frac{m+1}{m-k+1} \cdot \left(1 - \frac{\theta}{m+1}\right)^{k+t-1} \cdot \text{NegOcc}\left(t \middle| m, k, \frac{m\theta}{1-\theta+m}\right)$$

$$= \frac{m+1}{m-k+1} \cdot \left(1 - \frac{\theta}{m+1}\right)^{k+t} \cdot \text{NegOcc}\left(t \middle| m, k, \frac{m\theta}{1-\theta+m}\right),$$

which was to be shown. ■

**Proof of Theorem 8b:** Using the first recursive equation in Theorem 5 we obtain:

$$\text{NegOcc}(t|m, k+1, \theta) = \theta \cdot \frac{m-k}{m} \cdot \text{Occ}(k|k+t, m, \theta)$$

$$= \theta \cdot \frac{m-k}{m} \cdot \begin{bmatrix} \theta \cdot \frac{m-k+1}{m} \cdot \text{Occ}(k-1|k+t-1, m, \theta) \\ + \left(1 - \theta \cdot \frac{m-k}{m}\right) \cdot \text{Occ}(k|k+t-1, m, \theta) \end{bmatrix}$$

$$= \begin{bmatrix} \theta \cdot \frac{m-k}{m} \cdot \text{NegOcc}(t|m, k, \theta) \\ + \left(1 - \theta \cdot \frac{m-k}{m}\right) \cdot \text{NegOcc}(t-1|m, k+1, \theta) \end{bmatrix}.$$

Applying this recursive equation repeatedly to reduce the occupancy parameter in the second term gives us the desired result. (This result can be formally obtained by using induction with the present recursive equation.) ■

**Proof of Theorem 8c:** Using the differential equation in Theorem 5 we obtain:

$$\frac{\partial}{\partial\theta} \text{NegOcc}(t|m, k, \theta) = \frac{m-k+1}{m} \cdot \frac{\partial}{\partial\theta} [\theta \cdot \text{Occ}(k-1|k+t-1, m, \theta)]$$

$$= \frac{m-k+1}{m} \cdot \begin{bmatrix} \text{Occ}(k-1|k+t-1, m, \theta) \\ +\theta \cdot \frac{\partial}{\partial\theta} \text{Occ}(k-1|k+t-1, m, \theta) \end{bmatrix}$$

$$= \frac{m-k+1}{m} \cdot \begin{bmatrix} \text{Occ}(k-1|k+t-1, m, \theta) \\ -\theta \cdot \frac{m-k+1}{m} \cdot \text{Occ}(k-1|k+t-2, m, \theta) \\ +\theta \cdot \frac{m-k+2}{m} \cdot \text{Occ}(k-2|k+t-2, m, \theta) \end{bmatrix}$$

$$= \frac{m-k+1}{m} \cdot \begin{bmatrix} \text{Occ}(k-1|k+t-1, m, \theta) \\ -\text{NegOcc}(t-1|m,k,\theta) \\ +\text{NegOcc}(t|m,k-1,\theta) \end{bmatrix}$$

$$= \frac{1}{\theta} \cdot \text{NegOcc}(t|m,k,\theta) + \frac{m-k+1}{m} \cdot \begin{bmatrix} \text{NegOcc}(t|m,k-1,\theta) \\ -\text{NegOcc}(t-1|m,k,\theta) \end{bmatrix},$$

which was to be shown. ■

**Proof of Theorem 9a:** To facilitate our analysis, define the updated probability parameter $\theta' = m\theta/(1-\theta+m)$ which satisfies the equation $m \cdot (1-\theta')/\theta' = (m+1) \cdot (1-\theta)/\theta$. Also, recall the function $\Psi$ defined in the proof of Theorem 6b. Using the alternative form of the negative occupancy distribution, it is possible to show that:

$$\frac{\text{NegOcc}(r|m+1,k,\theta)}{\text{NegOcc}(r|m,k,\theta)} = \Psi(m,k) \cdot \frac{\Pi\left(k+t-1, k-1, m \cdot \frac{1-\theta'}{\theta'}\right)}{\Pi\left(k+t-1, k-1, m \cdot \frac{1-\theta}{\theta}\right)} \leq 1.$$

(We have equality if $k = 1$ and $\theta = 1$, and we have strict inequality if $k > 1$ or $\theta < 1$.) Thus, with the same conditions on strict inequality, we have:

$$F(t|m+1,k,\theta) = \sum_{r=0}^{t} \text{NegOcc}(r|m+1,k,\theta) \leq \sum_{r=0}^{t} \text{NegOcc}(r|m,k,\theta) = F(t|m+1,k,\theta).$$

Thus, for all $m \leq m'$ we have:

$$\begin{aligned} F(t|m',k,\theta) &= F(t|m,k,\theta) + \sum_{i=0}^{m'-m-1} \big(F(t|m+i+1,k,\theta) - F(t|m+i,k,\theta)\big) \\ &\leq F(t|m,k,\theta), \end{aligned}$$

where the inequality is strict if $m' > m$ and $k > 1$ or $\theta < 1$. ■

**Proof of Theorem 9b:** Applying the second recursive equation in Theorem 8 gives:

$$\begin{aligned} F(t|m,k+1,\theta) &= \sum_{r=0}^{t} \text{NegOcc}(r|m,k+1,\theta) \\ &= \sum_{r=0}^{t} \theta \cdot \frac{m-k}{m} \cdot \sum_{i=0}^{r} \left(1 - \theta \cdot \frac{m-k}{m}\right)^{i} \cdot \text{NegOcc}(r-i|m,k,\theta) \\ &= \theta \cdot \frac{m-k}{m} \sum_{i=0}^{t} \left(1 - \theta \cdot \frac{m-k}{m}\right)^{i} \sum_{r=i}^{t} \text{NegOcc}(r-i|m,k,\theta) \\ &= \theta \cdot \frac{m-k}{m} \sum_{i=0}^{t} \left(1 - \theta \cdot \frac{m-k}{m}\right)^{i} F(t-i|m,k,\theta) \end{aligned}$$

$$\leq \theta \cdot \frac{m-k}{m} \sum_{i=0}^{t} \left(1 - \theta \cdot \frac{m-k}{m}\right)^{i} F(t|m,k,\theta)$$

$$= \left[1 - \left(1 - \theta \cdot \frac{m-k}{m}\right)^{t+1}\right] \cdot F(t|m,k,\theta)$$

$$\leq F(t|m,k,\theta).$$

Thus, for all $k \leq k'$ we have:

$$F(t|m,k',\theta) = F(t|m,k,\theta) + \sum_{i=0}^{k'-k-1} \big(F(t|m,k+i+1,\theta) - F(t|m,k+i,\theta)\big)$$

$$\leq F(t|m,k,\theta) + \sum_{i=0}^{k'-k-1} \left[1 - \left(1 - \theta \cdot \frac{m-k-i}{m}\right)^{t+1}\right] \cdot F(t|m,k+i,\theta)$$

$$\leq F(t|m,k,\theta),$$

where the inequality is strict if $k < k'$. ■

**Proof of Theorem 9c:** The likelihood-ratio function for the negative occupancy distribution (comparing parameter values $\theta$ and $\theta'$) is:

$$\frac{\text{NegOcc}(t|m,k,\theta')}{\text{NegOcc}(t|m,k,\theta)} \overset{t}{\propto} \left(\frac{1-\theta'}{1-\theta}\right)^{t} \cdot \frac{\Pi\left(k+t-1, k-1, m \cdot \frac{1-\theta'}{\theta'}\right)}{\Pi\left(k+t-1, k-1, m \cdot \frac{1-\theta}{\theta}\right)}.$$

For $\theta' \geq \theta$ function is monotone non-increasing in $t$, and for $\theta' > \theta$ it is monotone increasing in $t$. The stochastic dominance result follows directly from this monotonicity property. ■

**Proof of Theorem 10:** Proof is analogous to Theorem 4(a), but we set it out explicitly to avoid any confusion. Applying Lemma 4 gives:

$$\lim_{\phi\to\infty} \text{Spillage}(r|n,k,\phi) = \lim_{\phi\to\infty} \binom{n}{k+r} \cdot \phi^{n-k-r} \cdot \frac{S(k+r,k)}{S(n,k,\phi)}$$

$$= S(k+r,k) \cdot \lim_{\phi\to\infty} \binom{n}{k+r} \cdot \frac{\phi^{n-k}}{S(n,k,\phi)} \times \lim_{\phi\to\infty} \frac{1}{\phi^{r}}$$

$$= S(k+r,k) \cdot \binom{n}{k+r} / \binom{n}{k} \times \mathbb{I}(r=0) = \mathbb{I}(r=0),$$

which was to be shown. ■

**Proof of Theorem 11:** Using the law of total probability we have:

$$p_K(k|m,\theta) = \sum_{n=0}^{\infty} p_N(n) \cdot \text{Occ}(k|n,m,\theta)$$

$$= \sum_{n=0}^{\infty} p_N(n) \cdot \binom{m}{k} \sum_{i=0}^{k} \binom{k}{i} (-1)^{k-i} \cdot \left(1 - \theta \cdot \frac{m-i}{m}\right)^n$$

$$= \binom{m}{k} \sum_{i=0}^{k} \binom{k}{i} (-1)^{k-i} \sum_{n=0}^{\infty} p_N(n) \cdot \left(1 - \theta \cdot \frac{m-i}{m}\right)^n$$

$$= \binom{m}{k} \sum_{i=0}^{k} \binom{k}{i} (-1)^{k-i} G_N \left(1 - \theta \cdot \frac{m-i}{m}\right),$$

which was to be shown. ■

**PROOF OF THEOREM 12:** This theorem is an application of Theorem 11 where $N \sim \text{Bin}(n, \theta)$. This distribution has probability generating function $G_N(z) = (1 - \theta + \theta z)^n$ so we have:

$$G_N \left(1 - \gamma \cdot \frac{m-i}{m}\right) = \left(1 - \theta + \theta \left(1 - \gamma \cdot \frac{m-i}{m}\right)\right)^n = \left(1 - \gamma\theta \cdot \frac{m-i}{m}\right)^n.$$

Applying Theorem 3 we have:

$$\mathbb{P}(K_N = k | m, \theta) = \binom{m}{k} \sum_{i=0}^{k} \binom{k}{i} (-1)^{k-i} G_N \left(1 - \gamma \cdot \frac{m-i}{m}\right)$$

$$= \binom{m}{k} \sum_{i=0}^{k} \binom{k}{i} (-1)^{k-i} \left(1 - \gamma\theta \cdot \frac{m-i}{m}\right)^n = \text{Occ}(k | n, m, \gamma\theta),$$

which establishes the main theorem. The special case holds trivially by substitution. ■

**PROOF OF THEOREM 13:** This result is an application of Theorem 11 where $N \sim \text{Pois}(\lambda)$. This distribution has probability generating function $G_N(z) = \exp(\lambda(z - 1))$ so we have:

$$G_N \left(1 - \theta \cdot \frac{m-i}{m}\right) = \exp\left(-\theta\lambda \cdot \frac{m-i}{m}\right).$$

Applying Theorem 3 and using the binomial theorem we have:

$$\mathbb{P}(K_N = k | m, \theta) = \binom{m}{k} \sum_{i=0}^{k} \binom{k}{i} (-1)^{k-i} G_N \left(1 - \theta \cdot \frac{m-i}{m}\right)$$

$$= \binom{m}{k} \sum_{i=0}^{k} \binom{k}{i} (-1)^{k-i} \exp\left(-\theta\lambda \cdot \frac{m-i}{m}\right)$$

$$= \exp(-\theta\lambda) \binom{m}{k} \sum_{i=0}^{k} \binom{k}{i} (-1)^{k-i} \exp(\theta\lambda/m)^i$$

$$= \exp(-\theta\lambda) \binom{m}{k} (\exp(\theta\lambda/m) - 1)^k$$

$$= \exp(-\theta\lambda/m)^m \binom{m}{k} (\exp(\theta\lambda/m) - 1)^k$$
$$= \binom{m}{k} (1 - \exp(-\theta\lambda/m))^k \exp(-\theta\lambda/m)^{m-k}$$
$$= \text{Bin}\left(k \middle| m, 1 - \exp\left(-\frac{\lambda\theta}{m}\right)\right),$$

which establishes the first equation in the theorem. The second equation follows immediately by substitution of the parameter value $\lambda = m \cdot |\ln(1-\phi)|/\theta$. ■

**PROOF OF THEOREM 14:** Let $\phi = m \cdot (1-\gamma)/\gamma$ and $\phi' = m \cdot (1-\gamma\theta)/\gamma\theta$ and note that:

$$\phi' - \phi = m \cdot \frac{1-\gamma\theta}{\gamma\theta} - m \cdot \frac{1-\gamma}{\gamma} = \frac{m}{\gamma} \cdot \frac{1-\theta}{\theta}.$$

By moving the noncentrality parameter using Lemma 3, we have:

$$S(k+t-1, k-1, \phi') = \sum_{r=0}^{t} \binom{k+t-1}{r} \left(\frac{m}{\gamma} \cdot \frac{1-\theta}{\theta}\right)^r \cdot S(k+t-r-1, k-1, \phi).$$

Hence, we have:

$$\text{NegOcc}(t|m,k,\gamma\theta) = \frac{(\gamma\theta)^{k+t}}{m^{k+t}} \cdot (m)_k \cdot S(k+t-1, k-1, \phi')$$
$$= \frac{(\gamma\theta)^{k+t}}{m^{k+t}} \cdot (m)_k \cdot \sum_{r=0}^{t} \binom{k+t-1}{r} \left(\frac{m}{\gamma} \cdot \frac{1-\theta}{\theta}\right)^r \cdot S(k+t-r-1, k-1, \phi)$$
$$= \sum_{r=0}^{t} \binom{k+t-1}{r} (1-\theta)^r \cdot \theta^{k+t-r} \cdot \frac{\gamma^{k+t-r}}{m^{k+t-r}} \cdot (m)_k \cdot S(k+t-r-1, k-1, \phi)$$
$$= \sum_{r=0}^{t} \text{NegBin}(r|k+t-r, 1-\theta) \cdot \text{NegOcc}(t-r|m,k,\gamma)$$
$$= \sum_{r=0}^{t} \text{NegBin}(t-r|k+r, 1-\theta) \cdot \text{NegOcc}(r|m,k,\gamma),$$

which establishes the main theorem. The special case holds trivially by substitution. ■

**PROOF OF THEOREM 15:** This theorem can be established as a consequence of the law of total probability, using the known distributions established in the body of the paper:

$$\text{Bin}(s|n,\theta) = \mathbb{P}(n_{\text{eff}} = s|n,m,\theta)$$
$$= \sum_{k=0}^{s} \mathbb{P}(n_{\text{eff}} = s, K = k|n,m,\theta)$$

$$= \sum_{k=0}^{s} \mathbb{P}(n_{\text{eff}} = s | K = k, n, m, \theta) \cdot \mathbb{P}(K = k | n, m, \theta)$$

$$= \sum_{k=0}^{s} \text{Spillage}\left(s - k \middle| n, k, m \cdot \frac{1-\theta}{\theta}\right) \cdot \text{Occ}(k | n, m, \theta)\,,$$

which was to be shown. ■